

\input amstex
\documentstyle{amsppt}

\input label.def
\input degt.def
\def\mnote#1{}



\input epsf
\def\picture#1{\epsffile{#1-bb.eps}}
\def\cpic#1{$\vcenter{\hbox{\picture{#1}}}$}

\def\ie{\emph{i.e.}}
\def\eg{\emph{e.g.}}
\def\cf.{\emph{cf\.}}
\def\via{\emph{via}}

\def\viz.{\emph{viz}\.}

{\catcode`\@11
\gdef\proclaimfont@{\sl}}

\newtheorem{Remark}\Remark\thm\endAmSdef
\conjecture\thm\endproclaim
\problem\thm\endAmSdef
\newhead\subsubsection\subsubsection\endsubhead
\def\paragraph{\subsubsection{}}

\def\dash{\item"\hfill--\hfill"}
\def\Dashes{\widestnumber\item{--}\roster}
\def\endDashes{\endroster}

\loadbold
\def\bA{\bold A}
\def\bD{\bold D}
\def\bE{\bold E}

\def\tA#1{\smash{\tilde\bA#1}}
\def\tD#1{\smash{\tilde\bD#1}}
\def\tE#1{\smash{\tilde\bE#1}}

\let\bc\beta
\let\MB=U

\let\splus\oplus

\def\FF{\Bbb F}

\def\CG#1{\Z_{#1}}
\def\BG#1{\Bbb B_{#1}}
\def\RBG#1{\bar\Bbb B_{#1}}
\def\DG#1{\Bbb D_{#1}}

\def\AG#1{\Bbb A_{#1}}
\def\TG#1{\Bbb T_{#1}}

\def\GL{\operatorname{\text{\sl GL\/}}}
\def\SL{\operatorname{\text{\sl SL\/}}}

\def\0{{\setbox0\hbox{0}\hbox to\wd0{\hss\rm--\hss}}}

\let\onto\twoheadrightarrow

\def\CC{\Cal C}

\def\B{\bar B}

\def\H{\bar H}
\def\F{\bar F}
\def\dR{\partial R}

\def\Cp#1{\Bbb P^{#1}}
\def\term#1-{$\DG{#1}$-}

\let\Ga\alpha
\let\Gb\beta

\let\Gs\sigma
\let\Gr\rho
\let\Gd\delta
\def\1{^{-1}}
\let\LL\Lambda

\def\ls|#1|{\mathopen|#1\mathclose|}
\let\<\langle
\let\>\rangle

\def\Aut{\operatorname{Aut}}

\def\ord{\operatorname{ord}}
\def\Sk#1{\operatorname{Sk}#1}

\def\tabstrut{\vrule height9.5pt depth2.5pt}
\def\exstrut{\omit\vrule height2pt\hss\vrule}

\def\GAP{{\tt GAP}}
\def\FRAG{fig.+e6}

\def\dobifrag#1{#1,\bar#1\bifragneg}
\def\dobifrag#1{#1}
\def\fragi(#1)#2{\ref{\FRAG}(#1)--$#2$}
\def\bifragi(#1)#2{\fragi(#1){\dobifrag{#2}}}
\def\fragii(#1)#2{\ref{fig.e6}(#1)--$#2$}
\def\bifragii(#1)#2{\fragii(#1){\dobifrag{#2}}}

\def\inserthyphen{\ifcat\next a-\fi\ignorespaces}
\let\BLACK\bullet
\let\WHITE\circ
\def\CROSS{\vcenter{\hbox{$\scriptstyle\mathord\times$}}}
\let\STAR*
\def\pblack-{$\BLACK$\futurelet\next\inserthyphen}
\def\pwhite-{$\WHITE$\futurelet\next\inserthyphen}
\def\pcross-{$\CROSS$\futurelet\next\inserthyphen}
\def\pstar-{$\STAR$\futurelet\next\inserthyphen}
\def\black{\protect\pblack}
\def\white{\protect\pwhite}
\def\cross{\protect\pcross}
\def\star{\protect\pstar}
\def\COUNT#1{\mathord\#_{#1}}
\def\nblack{\COUNT\BLACK}
\def\nwhite{\COUNT\WHITE}

\def\nstar{\COUNT\STAR}

\def\No{no\.~}
\def\Nos{nos\.~}

\def\notabelian{$^*$}
\def\notB{$^{**}$}
\def\notknown{$^{?}$}
\def\notknownB{$^{?\!{*}}$}
\let\MARK\relax
\def\defMARK#1{\gdef\MARK{\llap{#1}\global\let\MARK\relax}\MARK\ignorespaces}
\def\deftable{\def\1{\defMARK\notabelian}\def\2{\defMARK\notB}
 \def\3{\defMARK\notknownB}\def\4{\defMARK\notknown}}

\topmatter

\author
Alex Degtyarev
\endauthor


\title
Plane sextics with a type $\bold E_6$ singular point
\endtitle

\address
Department of Mathematics,
Bilkent University,
06800 Ankara, Turkey
\endaddress

\email
degt\@fen.bilkent.edu.tr
\endemail

\abstract
We give a classification up to equisingular deformation
and compute the fundamental
groups of maximizing plane sextics with a type~$\bE_6$
singular point.
\endabstract



\endtopmatter

\document

\section{Introduction}

\subsection{The subject}
This paper concludes the series~\cite{dessin.e7},
\cite{dessin.e8}, where we give a complete
deformation classification
and compute the fundamental
groups of maximizing irreducible plane sextics with an $\bE$ type
singular point. (With the common abuse of the language, by the
\emph{fundamental group} of a curve $B\subset\Cp2$ we mean the
group $\pi_1(\Cp2\sminus B)$ of its complement.)
Here, we consider sextics $B\subset\Cp2$
satisfying the following
conditions:
\roster
\item"$(*)$"
$B$ has simple (\ie, $\bA$--$\bD$--$\bE$) singularities only,
\item""
$B$ has a distinguished singular point~$P$ of type~$\bE_6$, and
\item""
$B$ has no singular points of type~$\bE_7$ or~$\bE_8$.
\endroster
(Singular points of type~$\bE_7$ and~$\bE_8$ are excluded in order
to reduce the lists.
Sextics with a type~$\bE_8$ point are considered
in~\cite{dessin.e8},
and irreducible sextics with a type~$\bE_7$ point are considered
in~\cite{dessin.e7}. Reducible sextics with a type~$\bE_7$ point,
as well as the more involved case of a distinguished $\bD$ type
point, may appear elsewhere.)

Recall that a plane sextic~$B$ with simple singularities only is
called \emph{maximizing} if the total Milnor number $\mu(B)$
assumes the maximal possible value~$19$. It is well known that
maximizing sextics are defined over algebraic number fields
(as they
are related to singular $K3$-surfaces).
Furthermore, such sextics are rigid: two maximizing sextics are
equisingular deformation equivalent if and only if they are
related by a projective transformation.

Another important class is formed by the so called sextics of
\emph{torus type}, \ie, those whose equation can be represented in
the form $f_2^3+f_3^2=0$, where $f_2$ and~$f_3$ are certain
homogeneous polynomials of degree~$2$ and~$3$, respectively. (This
property turns out to be equisingular deformation invariant.)
Each sextic~$B$ of torus type can be perturbed to Zariski's famous six
cuspidal sextic~\cite{Zariski.group}, which is obtained when~$f_2$
and~$f_3$ above are sufficiently generic. Hence, the group
$\pi_1(\Cp2\sminus B)$ factors to the \emph{reduced braid group}
$\bar\BG3:=\BG3/(\Gs_1\Gs_2)^3$; in particular, it is never
finite.
(The existence\mnote{a historical remark and reference to
\cite{Snyder} added}
of two distinct families of irreducible
six cuspidal sextics, those of and not of torus type, was first
stated by Del Pezzo and then proved by B.~Segre, see
\eg~\cite{Snyder, page~407}. Later Zariski~\cite{Zariski.group}
showed that the two families differ by the fundamental groups.)

A representation of the equation of~$B$ in the form
$f_2^3+f_3^2=0$ is called a \emph{torus structure}. The points of
intersection of the conic $\{f_2=0\}$ and cubic $\{f_3=0\}$
are always singular for~$B$; they are called \emph{inner} (with
respect to the given torus structure), whereas
the other singular points are called \emph{outer}. In the listing
below, we indicate sextics of torus type by representing their
sets of singularities in the form
$$
(\text{inner singularities})\splus\text{outer singularities}.
$$
An exception is the set of singularities
$\bE_6\splus\bA_5\splus4\bA_2$, which is always of torus type and
admits four distinct torus
structures.

Formally, the deformation classification of plane sextics with
simple singularities can be reduced to a purely arithmetical
problem, see~\cite{JAG}, and for maximizing sextics this latter
problem has been completely solved, see~\cite{Shimada}
and~\cite{Yang}, in the sense that all deformation classes have
been enumerated. Unfortunately, this approach, based on the theory
of $K3$-surfaces and the global Torelli theorem, is not
constructive and very little is known about the geometry of the
curves. (Sporadic examples using explicit equations are scattered
in the literature.) Here, we use another approach, suggested
in~\cite{degt.kplets} and~\cite{dessin.e7}: a plane curve~$B$ with
a sufficiently deep (with respect to the degree) singularity is
reduced to a trigonal curve~$\B$ in an appropriate Hirzebruch
surface. If $B$ is a maximizing sextic (with a triple point), then
$\B$ is a maximal trigonal curve; hence it can be studied using
Grothendieck's \emph{dessin d'enfants} of its functional
$j$-invariant. At the end, we obtain an explicit geometric
description of~$\B$ and~$B$ (rather than equation);
among other things, this description is suitable for computing the
braid monodromy and hence the fundamental group of the curves.

\subsection{Results}
The principal results of the present paper are
Theorems~\ref{th.classification} (the classification)
and~\ref{th.group}, \ref{th.pert.nontorus}
(the computation of the fundamental group) below.

\theorem\label{th.classification}
Up to projective transformation \rom(equivalently, up to
equisingular deformation\rom)
there are $93$
maximizing plane
sextics satisfying condition~$(*)$ above, realizing $71$
combinatorial sets of singularities\rom; of them, $53$ sextics
\rom($40$ sets of singularities\rom) are irreducible, see
Table~\ref{tab.e6} on page~\rom{\pageref{tab.e6}}, and $40$ sextics
\rom($32$ sets of singularities\rom) are reducible, see
Table~\ref{tab.e6-r} on page~\rom{\pageref{tab.e6-r}}.
\endtheorem

In Theorem~\ref{th.classification},
one set of singularities is common:
$\bE_6\splus\bA_9\splus\bA_4$ is realized by
three irreducible and one reducible sextics. This theorem is
proved in Section~\ref{S.classification}; for more details, see
comments to the tables in Subsection~\ref{s.list}.

Among the irreducible sextics in Theorem~\ref{th.classification},
twelve
(eight sets of singularities) are of torus type.
Seven of them (four sets of singularities)
are `new' in the sense that they have not been
extensively studied before.

\theorem\label{th.group}
Let $B\subset\Cp2$ be
an irreducible\mnote{`irreducible' inserted}
maximizing sextic that
satisfies condition~$(*)$
above and is not of torus type. If the
set of singularities of~$B$ is
$2\bE_6\splus\bA_4\splus\bA_3$ \rom(\Nos$4$ and~$5$ in
Table~\ref{tab.e6}\rom), then
$\pi_1(\Cp2\sminus B)=\SL(2,\FF_5)\rtimes\CG6$
\rom(see~\eqref{rep.hloop} or~\eqref{rep.h2loop} for the
presentations\rom). Otherwise, $\pi_1(\Cp2\sminus B)=\CG6$.
\endtheorem

\theorem\label{th.pert.nontorus}
Let~$B$ be a sextic as in Theorem~\ref{th.classification}, and
let~$B'$ be a proper irreducible perturbation of~$B$ that is not
of torus type. Then $\pi_1(\Cp2\sminus B')=\CG6$.
\endtheorem

Theorems~\ref{th.group} and~\ref{th.pert.nontorus} are proved in
Section~\ref{S.computation} and
Subsection~\ref{pf.pert.nontorus}, respectively.
The theorems substantiate my conjecture asserting that the fundamental
group of an irreducible plane sextic that has simple singularities
only and
is not of torus type is
finite. Recall\mnote{motivation and current state of the
conjecture discussed; `simple singularities only' inserted}
that, originally, the conjecture was motivated by a certain
experimental evidence (which has now been extended) and the fact
that the abelianization of the fundamental group of an irreducible
sextic that is not of torus type is finite (which is
a restatement of
the proved part of
the so called Oka conjecture, see~\cite{degt.Oka}). At present,
the conjecture is essentially settled for sextics with a triple
singular point (the case of a $\bD$ type point is considered in a
forthcoming paper).

In Section~\ref{S.computation}, we also write down presentations
of the fundamental groups of all other sextics as in
Theorem~\ref{th.classification}, and in Section~\ref{S.pert}, we
consider their perturbations.
In particular,\mnote{\vskip-.7cm
the rest of introduction edited as there has
been a new development; new Theorem \ref{th.torus}
inserted, Problem 1.2.4
removed, and Conjecture 1.2.6 converted to Theorem~\ref{th.pert.torus}}
we prove the following theorem, computing the groups of all `new'
irreducible sextics of torus type.

\theorem\label{th.torus}
The\mnote{new theorem}
fundamental group of a sextic with the set of singularities
$(\bE_6\splus\bA_{11})\splus\bA_2$,
$(\bE_6\splus\bA_8\splus\bA_2)\splus\bA_3$, or
$(\bE_6\splus\bA_8\splus\bA_2)\splus\bA_2\splus\bA_1$
\rom(\Nos$12$, $13$, $18$, and~$40$ in Table~\ref{tab.e6}\rom)
is isomorphic to the
reduced braid group $\bar\BG3$. The groups of sextics with the set
of singularities
$(\bE_6\splus\bA_5\splus2\bA_2)\splus\bA_4$
\rom(\Nos$9$ and~$41$\rom) are \emph{not} isomorphic to
$\bar\BG3$\rom; their presentations are found
in~\ref{s.e6+a5+2a2+a4} and~\iref{s.group.hexagon}{41},
respectively.
\endtheorem


\problem
Are the fundamental groups of the
two
sextics with the set of singularities
$(\bE_6\splus\bA_5\splus2\bA_2)\splus\bA_4$ (\Nos$9$ and~$41$ in
Table~\ref{tab.e6}) isomorphic to each other?
(A similar question still stands for the
sextics with the set of singularities
$(2\bE_6\splus\bA_5)\splus\bA_2$, \Nos$6$ and~$7$
in Table~\ref{tab.e6},
see~\cite{degt.e6} and~\cite{OkaPho}.)
\endproblem

\theorem\label{th.pert.torus}
Let~$B$\mnote{converted to a theorem}
be a sextic as in Theorem~\ref{th.classification}, and
let~$B'$ be a proper irreducible perturbation of~$B$ that is
of torus type. Then, with the following few exceptions\rom:
\Dashes
\dash
a perturbation
$\bE_6\splus\bA_5\splus4\bA_2\to
\text{\rm(a curve of weight~$8$)}$,
see~\cite{degt.8a2}, and
\dash
a perturbation that can further be perturbed to a curve~$B''$ with
the set of singularities $(6\bA_2)\splus4\bA_1$,
see~\cite{degt.e6} and~\cite{degt.8a2},
\endDashes
the group $\pi_1(\Cp2\sminus B')$ is the reduced braid group
$\bar\BG3$.
\endtheorem

(Here, the \emph{weight} of a sextic is understood in the sense
of~\cite{degt.Oka},\mnote{weight defined}
as the total weight if all its singular points, where the weight
$w(P)$ of a singular point~$P$
is defined \via\ $w(\bA_{3k-1})=k$, $w(\bE_6)=2$, and
$w(P)=0$ otherwise.
The
fundamental group of a sextic of
weight~$\ge8$ is much larger than~$\bar\BG3$, as it has a larger
Alexander polynomial.
In the second exceptional case in the statement, one has
$\pi_1(\Cp2\sminus B')=\BG4/\Gs_1^2\Gs_2\Gs_3^2\Gs_2$,
see~\cite{degt.e6} and~\cite{degt.8a2}.)
Theorem~\ref{th.pert.torus}
is proved in Subsection~\ref{pf.pert.torus}.

\subsection{Acknowledgements}
I am grateful to E.~Artal Bartolo,\mnote{acknowledgement added}
who helped me to identify some
of the groups of curves of torus type,
thus making the statements more
complete,
and to I.~Dolgachev for his enlightening remarks concerning
history of six cuspidal sextics.

\section{The classification\label{S.classification}}

\subsection{The settings}\label{s.settings}
We remind briefly some of the results of~\cite{dessin.e7},
concerning the
construction and the classification of plane sextics satisfying~$(*)$.
For details on maximal trigonal curves and their skeletons,
see~\cite{degt.kplets} or~\cite{dessin.e7}. We denote
by~$\Sigma_k$, $k>0$, the geometrically ruled rational surface
with an exceptional section~$E$ of self-intersection~$-k$.

\proposition\label{1-1.e6}
There is a natural
bijection~$\phi$, invariant under equisingular deformations,
between Zariski open \rom(in each equisingular
stratum\rom) subsets of the following two sets\rom:
\roster
\item\local{e6.sextic}
plane sextics~$B$ with a distinguished
type~$\bE_6$ singular point~$P$, and
\item
trigonal curves $\B\subset\Sigma_4$ with a distinguished
type~$\tA{_5}$ singular fiber $\F$.
\endroster
A sextic~$B$ is irreducible if and only if so is $\B=\phi(B)$, and
$B$ is maximizing if and only if $\B$ is
maximal and stable, \ie, has no singular fibers of
types~$\tA{_0^{**}}$, $\tA{_1^*}$, or~$\tA{_2^*}$ \rom(these
fibers are called \emph{unstable}\rom).
\endproposition

Up to fiberwise equisingular deformation (equivalently, up to
automorphism of~$\Sigma_k$), maximal trigonal curves
$\B\subset\Sigma_k$ are classified by their
\emph{skeletons} and \emph{type specifications}. The skeleton
$\Sk=\Sk_{\B}\subset S^2$
(which is defined as Grothendieck's \emph{dessin d'enfants} of the
functional $j$-invariant of~$\B$)
is an embedded connected bipartite graph
with all \black-vertices of valency~$\le3$\mnote{$\le3$ retained;
the fact that \emph{in this paper} the actual values are $3$ or
$1$ (since there are no type $\bE_8$ points) is stated right
before \eqref{eq.count}}
and all \white-vertices
of valency~$\le2$. The \black-vertices of valency~$\le2$ and
\white-vertices of valency~$1$ are called \emph{singular}; they
correspond to the unstable and type~$\tE{}$ singular fibers
of~$\B$. Besides, each $n$-gonal \emph{region} of~$\Sk$ (\ie,
connected component of the complement $S^2\sminus\Sk$) contains a
single singular fiber of~$\B$, which is of type~$\tA{_{n-1}}$
($\tA{_0^*}$ if $n=1$) or~$\tD{_{n+4}}$. The type specification is
the function choosing, for each singular vertex and each region
of~$\Sk$, whether the corresponding fiber is of type~$\tA{}$
or~$\tD{}$, $\tE{}$.

The skeleton and the type specification of a maximal
curve~$\B\subset\Sigma_k$ are subject to the relation
$$
\nblack+\nwhite(1)+\nblack(2)=2(k-t),
$$
where $t$ is the number of triple singular points of the curve,
$\nstar(n)$ is the number of \star-vertices of valency~$n$,
$*=\bullet$ or~$\circ$, and $\nstar$ is the total number of
\star-vertices. Any pair satisfying this relation gives rise to a
unique curve.

Under the assumptions of this paper ($\B$ has no unstable fibers
or fibers of type~$\tE{_7}$ or~$\tE{_8}$), all \white-vertices
of~$\Sk$ are of valency~$2$ and all its \black-vertices are of
valency~$3$ or~$1$, the latter corresponding to the type~$\tE{_6}$
singular fibers of~$\B$. Hence the vertex count above can be
simplified to
$$
\nblack=2(k-t).
\eqtag\label{eq.count}
$$
Furthermore, the \white-vertices can be ignored, with the
convention that a \white-vertex is to be understood at the middle
of each edge connecting two \black-vertices.

To summarize, proof of Theorem~\ref{th.classification} reduces to
the enumeration of all pairs
$(\Sk,\text{type specification})$, where
$\Sk\subset S^2$ is a connected graph with all vertices of
valency~$3$ or~$1$ and with a distinguished hexagonal region.

\subsection{The case of two type~$\bE_6$ points}
The only maximizing sextic with three type~$\bE_6$ singular points
(\No$1$ in Table~\ref{tab.e6}) is well known, see~\cite{degt.e6},
\cite{OkaPho}.
Assume that $B$ has two type~$\bE_6$ singular points. Then~$\Sk$
has one monovalent \black-vertex and, in view of~\eqref{eq.count},
it can be obtained by attaching the fragment
${\bullet}{\joinrel\relbar\joinrel\relbar\joinrel}{\bullet}$ at
the center of an edge of a regular $3$-graph~$\Sk'$
with two or four vertices, see~\cite{symmetric}.
All possibilities resulting in a
skeleton~$\Sk$ with a hexagonal region are listed in
Figure~\ref{fig.+e6}, where $\Sk'$ is shown in black and the
possible position of the insertion, in grey. The sextics obtained
are \Nos$2$--$8$ in Table~\ref{tab.e6}; all curves are irreducible
due to the existence of a type~$\tE{_6}$ singular fiber.

\midinsert
\centerline{\vbox{\halign{\hss#\hss&&\qquad\quad\hss#\hss\cr
\cpic{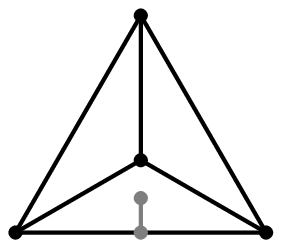}&
\cpic{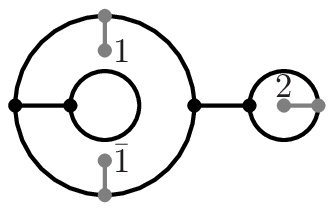}&
\cpic{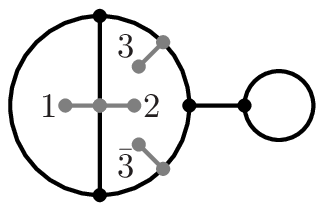}\cr
\noalign{\medskip}
(a)&(b)&(c)\cr
\crcr}}}
\figure
Two type~$\bE_6$ points
\endfigure\label{fig.+e6}
\endinsert

\Remark\label{rem.6-7}
For insertion~$2$ in Figure~\ref{fig.+e6}(c), the
skeleton~$\Sk$ has two hexagonal regions, resulting in two
sextics (\Nos$6$ and~$7$ in Table~\ref{tab.e6}). There are indeed
two distinct deformation families of sextics of torus type
with the set of singularities $(2\bE_6\splus\bA_5)\splus\bA_2$,
see~\cite{degt.e6}, \cite{OkaPho} and Remark~\ref{rem.4,5} below.
\endRemark

\subsection{The case of a hexagon with a loop}\label{s.hloop}
Till the rest of this section, assume that $P$ is the only
type~$\bE_6$ singular point of~$B$. Then $\Sk$ is a regular
$3$-graph with a distinguished hexagonal region~$\H$.
Such skeletons can be enumerated using~\cite{Beukers}; however, we
choose a more constructive descriptive approach.

Combinatorially, there are three possibilities for~$\H$:
\Dashes
\dash
hexagon with a loop, see Figure~\ref{fig.hloop}, left,
\dash
hexagon with a double loop, see Figure~\ref{fig.h2loop}, or
\dash
genuine hexagon, see Subsection~\ref{s.genuine} below.
\endDashes

\midinsert
\centerline{\cpic{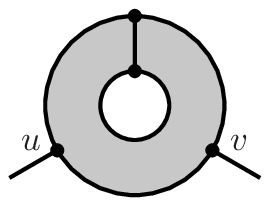}\qquad\qquad\cpic{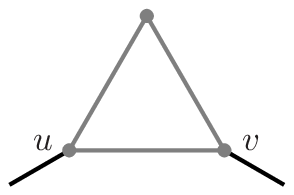}}
\figure
A hexagon with a loop
\endfigure\label{fig.hloop}
\endinsert

Assume that $\H$ is a hexagon with a loop, see the shaded area in
Figure~\ref{fig.hloop}, left. Removing a neighborhood of~$\H$
from~$\Sk$ and patching vertices~$u$, $v$ in the figure to a
single edge, one obtains another regular $3$-graph~$\Sk'$ with at
most four vertices, see~\cite{symmetric}. Conversely, $\Sk$ can be
obtained from~$\Sk'$ by inserting a fragment as in
Figure~\ref{fig.hloop}, left, at the middle of an edge of~$\Sk'$.
The essentially distinct possibilities for the position of the
insertion are shown in Figures~\ref{fig.e6} (irreducible curves)
and~\ref{fig.e6-r} (reducible curves; a reducibility criterion is
found in~\cite{degt.kplets}). To simplify the drawings, we
represent the insertion by a grey triangle, as in
Figure~\ref{fig.hloop}, right.

The resulting sextics are \Nos$9$--$32$ in Table~\ref{tab.e6} and
\Nos$1'$--$19'$ in Table~\ref{tab.e6-r}.

\midinsert
\centerline{\vbox{\halign{\hss#\hss&&\qquad\quad\hss#\hss\cr
\cpic{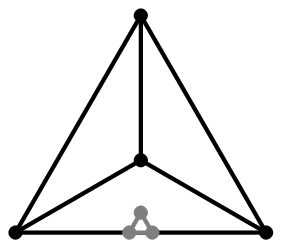}&\cpic{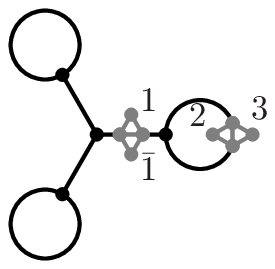}&\cpic{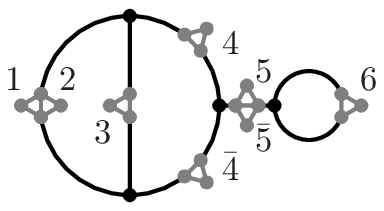}\cr
\noalign{\medskip}
(a)&(b)&(c)\cr
\crcr}}}
\bigskip
\centerline{\vbox{\halign{\hss#\hss&&\qquad\hss#\hss\cr
\cpic{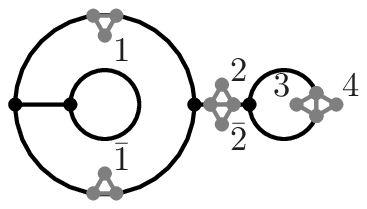}&\cpic{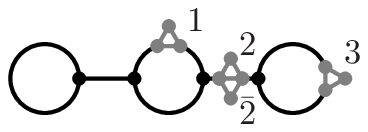}&\cpic{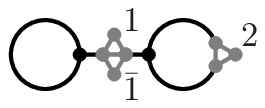}\cr
\noalign{\medskip}
(d)&(e)&(f)\cr
\crcr}}}
\figure
A hexagon with a loop: irreducible curves
\endfigure\label{fig.e6}
\endinsert

\midinsert
\centerline{\vbox{\halign{\hss#\hss&&\qquad\qquad\hss#\hss\cr
\cpic{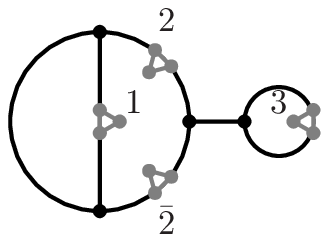}&\cpic{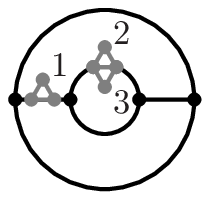}\cr
\noalign{\medskip}
(a)&(b)\cr
\crcr}}}
\bigskip
\centerline{\vbox{\halign{\hss#\hss&&\kern1.5em\hss#\hss\cr
\picture{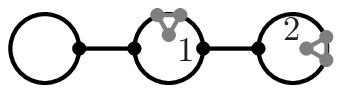}&\picture{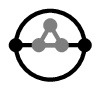}&\picture{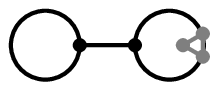}&\picture{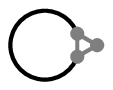}\cr
\noalign{\medskip}
(c)&(d)&(e)&(f)\cr
\crcr}}}
\figure
A hexagon with a loop: reducible curves
\endfigure\label{fig.e6-r}
\endinsert

\Remark\label{rem.27-28}
The curves in
pairs \Nos$27$, $28$, \Nos$31$, $32$, and \Nos$15'$, $16'$ in the
tables differ
by their type specifications:
the type~$\tD{_5}$ fiber can be chosen either inside one of the
`free' loops of~$\Sk'$
shown in the figure or inside the inner loop of the
hexagon. We assume that the latter possibility corresponds to
curves \Nos$28$, $32$ in Table~\ref{tab.e6} and \No$16'$ in
Table~\ref{tab.e6-r}.
\endRemark

\subsection{The case of a hexagon with a double loop}\label{s.h2loop}
Now, assume that the distinguished hexagon~$\H$ looks like the
outer region in Figure~\ref{fig.h2loop}, left. Each of the
remaining fragments~$A$, $B$ of~$\Sk$ has an odd number of
vertices, and the total number of remaining vertices is at
most four. Hence, one can assume that $A$ has one vertex and $B$
has at most three vertices. Then $A$ is a single loop and the graph can be
redrawn as shown in Figure~\ref{fig.h2loop}, right, where $\H$ is
represented by the shaded area. In other words, $\Sk$ can be
obtained from a regular $3$-graph~$\Sk'$ with two or four vertices
and with a loop, see~\cite{symmetric},
by replacing a loop with the fragment shown in
Figure~\ref{fig.h2loop}, right. The five possibilities
are listed in Figure~\ref{fig.e6-h}; the resulting sextics are
\Nos$33$--$38$ in Table~\ref{tab.e6}. (Using~\cite{degt.kplets},
one can easily show that the existence of a fragment as in
Figure~\ref{fig.h2loop}, right implies that the curve is
irreducible.)

\midinsert
\centerline{\cpic{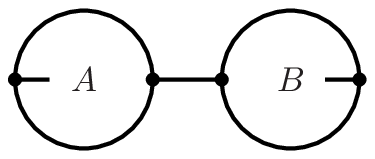}\qquad\qquad\cpic{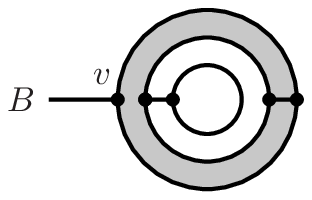}}
\figure
A hexagon with a double loop
\endfigure\label{fig.h2loop}
\endinsert

\midinsert
\centerline{\vbox{\halign{\hss#\hss&&\qquad\quad\hss#\hss\cr
\cpic{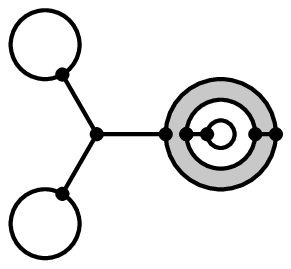}&\cpic{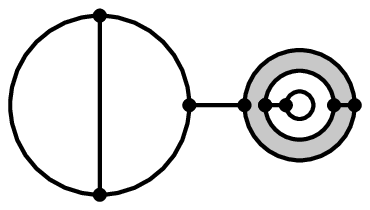}\cr
\noalign{\medskip}
(b)&(c)\cr
\crcr}}}
\bigskip
\centerline{\vbox{\halign{\hss#\hss&&\qquad\hss#\hss\cr
\cpic{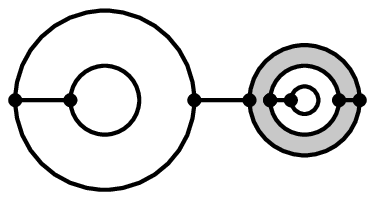}&\cpic{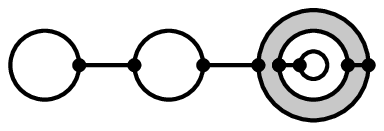}&\cpic{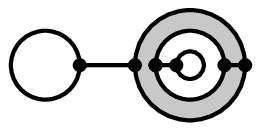}\cr
\noalign{\medskip}
(d)&(e)&(f)\cr
\crcr}}}
\figure
A hexagon with a double loop: the five skeletons
\endfigure\label{fig.e6-h}
\endinsert

\Remark\label{rem.37,38}
The skeleton in Figure~\ref{fig.e6-h}(f) has a symmetry
interchanging its two monogons and two pentagons. For this reason,
unlike the case described in Remark~\ref{rem.27-28},
items \Nos$37$ and~$38$ are realized by one deformation family
each.
\endRemark

\subsection{The case of a genuine hexagon}\label{s.genuine}
Finally, assume that $\H$ is a genuine hexagon, \ie, all six
vertices in the boundary~$\partial\H$ are pairwise distinct.
In other words, $\partial\H$ is the equator in~$S^2$, and $\Sk$
is obtained from~$\partial\H$ by completing it to a regular
$3$-graph by inserting at most two vertices and connecting edges
into one of the two hemispheres. All possibilities are listed in
Figures~\ref{fig.hex} (irreducible curves) and~\ref{fig.hex-r}
(reducible curves; a reducibility criterion is found
in~\cite{degt.kplets}); the resulting sextics are \Nos$39$--$42$
in Table~\ref{tab.e6} and \Nos$20'$--$38'$ in
Table~\ref{tab.e6-r}.

\midinsert
\def\pic(#1)#2{\vtop{\halign{##\cr
 \cpic{#2}\cr\noalign{\medskip}\hss(#1)\hss\cr}}}
\centerline{\pic(a){e6+a5+4a2}\qquad
 \pic(b){e6+a8+2a2+a1}\qquad
 \pic(c){e6+a5+2a2+a4}\qquad
 \pic(d){e6+a6+a4+a2+a1}}
\figure
A genuine hexagon: irreducible curves
\endfigure\label{fig.hex}
\endinsert

\midinsert
\def\pic(#1)#2{\vtop{\halign{##\cr
 \cpic{#2}\cr\noalign{\medskip}\hss\strut(#1)\hss\cr}}}
\centerline{\pic(a){e6+a5+2a2+a3+a1}\qquad
 \pic(b){e6+a4+2a3+a2+a1}\qquad
 \pic(c){e6+a7+a3+3a1}\qquad
 \pic(d){e6+a5+2a3+2a1}}
\bigskip
\centerline{\pic(e){e6+2a5+3a1}\qquad
 \pic(f){e6+a9+a3+a1}\qquad
 \pic(g){e6+a11+2a1-1}\qquad
 \pic(h){e6+a7+a5+a1}}
\bigskip
\centerline{\pic(i){e6+a11+2a1-2}\qquad
 \pic(j){e6+2a5+a3}\qquad
 \pic(k){e6+a5+a4+a3+a1}\qquad
 \pic(l){e6+a6+a5+2a1}}
\bigskip
\centerline{\pic(m){e6+a7+a4+2a1}\qquad
 \pic(n){e6+a7+a3+a2+a1}\qquad
 \pic(o){e6+a9+a2+2a1}}
\bigskip
\centerline{\pic(p){e6+d-1}\qquad
 \pic(q){e6+d-2}}
\figure
A genuine hexagon: reducible curves
\endfigure\label{fig.hex-r}
\endinsert

\subsection{The list}\label{s.list}
The results of the classification are collected in
Tables~\ref{tab.e6} (irreducible curves)
and~\ref{tab.e6-r} (reducible curves), where we list the
combinatorial types of singularities and references to the
corresponding figures. (For sextics of torus type, the inner
singularities are also indicated. Exception is \No$39$ in
Table~\ref{tab.e6}, which admits four distinct torus structures.)
Equal superscripts precede combinatorial types shared by several
items in the tables. (One set of singularities,
$\bE_6\splus\bA_9\splus\bA_4$ marked with~$^{10}$, is common for
both tables.)
The `Count' column lists the numbers
$(n_r,n_c)$ of real curves and pairs of complex conjugate curves.
The last two columns refer to the computation of the fundamental
group and indicate the parameters used (explained below).
Besides, the curves with nonabelian fundamental group
$\pi_1:=\pi_1(\Cp2\sminus B)$ are prefixed
with one of the following symbols:
\Dashes
\dash
\notabelian: the group $\pi_1$ is not abelian;
\dash
\notknown: the group~$\pi_1$ is not known to be abelian;
\dash
\notB: for curves of torus type,
$\pi_1\ne\bar\BG3$.\mnote{`\notknownB' option removed,
references in Table~\ref{tab.e6} edited}
\endDashes

\midinsert
\table\label{tab.e6}
Maximal sets of singularities with a type~$\bE_6$ point
represented by irreducible sextics
\endtable
\deftable

\def\no{}
\def\tref#1{\,\text{\ref{#1}}\,}
\def\tiref#1#2{\,\text{\iref{#1}{#2}}\,}
\def\tcite#1{\text{see \cite{#1}}}
\def\same{\afterassignment\dosame\count0=}
\def\dosame{\llap{$^{\the\count0\,}$}}
\centerline{\vbox{\offinterlineskip\halign{%
\tabstrut\ \ \ \hss#\hss\ \vrule&
\quad\,$#$\hss\quad\vrule&
\ \,\null#\hss\ \vrule&&\ \hss$#$\hss\ \vrule\cr
\noalign{\hrule}
\exstrut&&&&&\cr
\#\ \ &\text{Set of singularities}&
 \hss Figure&\text{Count}&\pi_1&\text{Parameters}\cr
\exstrut&&&&&\cr
\noalign{\hrule}
\exstrut&&&&&\cr
\2 1&(3\bE_6)\splus\bA_1&
 &(1,0)
 &\tcite{degt.e6}&\cr
\exstrut&&&&&\cr
\noalign{\hrule}
\exstrut&&&&&\cr
\2 2&(2\bE_6\splus2\bA_2)\splus\bA_3&
 \ref{fig.+e6}(a)&(1,0)
 &\tcite{degt.e6}&\cr
3&2\bE_6\splus\bA_7&
 \bifragi(b)1&(0,1)
 &\tref{s.group.loop}&(\0,\0,1,\0)\cr
\1 4&\same1 2\bE_6\splus\bA_4\splus\bA_3&
 \fragi(b)2&(1,0)
 &\tref{s.group.2loop}&l=4\cr
\1 5&\same1 2\bE_6\splus\bA_4\splus\bA_3&
 \fragi(c)1&(1,0)
 &\tref{s.group.loop}&(4,5,\0,\0)\cr
\2 6&\same2 (2\bE_6\splus\bA_5)\splus\bA_2&
 \fragi(c)2&(1,0)
 &\tcite{degt.e6}&\cr
\2 7&\same2 (2\bE_6\splus\bA_5)\splus\bA_2&
 \fragi(c)2&(1,0)
 &\tcite{degt.e6}&\cr
8&2\bE_6\splus\bA_6\splus\bA_1&
 \bifragi(c)3&(0,1)
 &\tref{s.group.2e6}&\cr
\exstrut&&&&&\cr
\noalign{\hrule}
\exstrut&&&&&\cr
\2 9&\same3 (\bE_6\splus\bA_5\splus2\bA_2)\splus\bA_4&
 \ref{fig.e6}(a)&(1,0)
 &\tref{s.e6+a5+2a2+a4}\!\!&(6,5,3,\0)\cr
10&\bE_6\splus\bA_{13}&
 \bifragii(b)1&(0,1)
 &\tref{s.group.loop}&(\0,\0,1,\0)\cr
11&\same4 \bE_6\splus\bA_{10}\splus\bA_3&
 \fragii(b)2&(1,0)
 &\tref{s.e6+a10+a3}&\cr
\1 12&(\bE_6\splus\bA_{11})\splus\bA_2&
 \fragii(b)3&(1,0)
 &\tref{s.e6+a11+a2}&\cr
\1 13&(\bE_6\splus\bA_8\splus\bA_2)\splus\bA_3&
 \fragii(c)1&(1,0)
 &\tref{s.loop.torus}\!\!&(9,4,3,\0)\cr
14&\same5 \bE_6\splus\bA_7\splus\bA_4\splus\bA_2&
 \fragii(c)2&(1,0)
 &\tref{s.group.loop}&(5,8,3,\0)\cr
15&\same6 \bE_6\splus\bA_5\splus2\bA_4&
 \fragii(c)3&(1,0)
 &\tref{s.group.loop}&(5,5,6,\0)\cr
16&\same7 \bE_6\splus\bA_8\splus\bA_4\splus\bA_1&
 \bifragii(c)4&(0,1)
 &\tref{s.group.loop}&(\0,\0,\0,1)\cr
17&\same8 \bE_6\splus\bA_{10}\splus\bA_2\splus\bA_1&
 \bifragii(c)5&(0,1)
 &\tref{s.group.loop}&(\0,\0,1,\0)\cr
\1 18&\same9 (\bE_6\splus\bA_8\splus\bA_2)\splus\bA_2\splus\bA_1\!\!\!&
 \fragii(c)6&(1,0)
 &\tref{s.loop.torus}\!\!&(9,3,\0,3)\cr
19&\bE_6\splus\bA_7\splus\bA_6&
 \bifragii(d)1&(0,1)
 &\tref{s.group.loop}&(\0,\0,\0,1)\cr
20&\same10 \bE_6\splus\bA_9\splus\bA_4&
 \bifragii(d)2&(0,1)
 &\tref{s.group.loop}&(\0,\0,1,\0)\cr
21&\bE_6\splus\bA_6\splus\bA_4\splus\bA_3&
 \fragii(d)3&(1,0)
 &\tref{s.group.loop}&(4,7,\0,\0)\cr
22&\same5 \bE_6\splus\bA_7\splus\bA_4\splus\bA_2&
 \fragii(d)4&(1,0)
 &\tref{s.group.loop}&(3,8,\0,5)\cr
23&\same4 \bE_6\splus\bA_{10}\splus\bA_3&
 \fragii(e)1&(1,0)
 &\tref{s.group.loop}&(\0,\0,\0,1)\cr
24&\bE_6\splus\bA_{12}\splus\bA_1&
 \bifragii(e)2&(0,1)
 &\tref{s.group.loop}&(\0,\0,1,\0)\cr
25&\same8 \bE_6\splus\bA_{10}\splus\bA_2\splus\bA_1&
 \fragii(e)3&(1,0)
 &\tref{s.group.loop}&(11,3,\0,2)\cr
26&\bE_6\splus\bD_{13}&
 \fragii(f)1&(1,0)
 &\tref{s.group.loop}&(\0,\0,1,\0)\cr
27&\same11 \bE_6\splus\bD_5\splus\bA_8&
 \bifragii(f)1&(0,1)
 &\tref{s.group.loop}&(\0,\0,1,\0)\cr
28&\same11 \bE_6\splus\bD_5\splus\bA_8&
 \bifragii(f)1&(1,0)
 &\tref{s.e6+d5+a8}&\cr
29&\bE_6\splus\bD_{11}\splus\bA_2&
 \fragii(f)2&(1,0)
 &\tref{s.group.loop}&(\0,\0,\0,1)\cr
30&\bE_6\splus\bD_7\splus\bA_6&
 \fragii(f)2&(1,0)
 &\tref{s.group.loop}&(\0,\0,\0,1)\cr
31&\same12 \bE_6\splus\bD_5\splus\bA_6\splus\bA_2&
 \fragii(f)2&(1,0)
 &\tref{s.group.loop}&(7,3,\0,\0)\cr
32&\same12 \bE_6\splus\bD_5\splus\bA_6\splus\bA_2&
 \fragii(f)2&(1,0)
 &\tref{s.e6+d5+a6+a2}&\cr
\exstrut&&&&&\cr
\noalign{\hrule}
\exstrut&&&&&\cr
33&\same10 \bE_6\splus\bA_9\splus\bA_4&
 \ref{fig.e6-h}(b)&(1,0)
 &\tref{s.group.2loop}&l=10\cr
34&\same13 \bE_6\splus\bA_6\splus\bA_4\splus\bA_2\splus\bA_1&
 \ref{fig.e6-h}(c)&(1,0)
 &\tref{s.group.2loop}&l=7\cr
35&\same6 \bE_6\splus\bA_5\splus2\bA_4&
 \ref{fig.e6-h}(d)&(1,0)
 &\tref{s.group.2loop}&l=6\cr
36&\same7 \bE_6\splus\bA_8\splus\bA_4\splus\bA_1&
 \ref{fig.e6-h}(e)&(1,0)
 &\tref{s.group.2loop}&l=9\cr
37&\bE_6\splus\bD_9\splus\bA_4&
 \ref{fig.e6-h}(f)&(1,0)
 &\tref{s.group.2loop}&\cr
38&\bE_6\splus\bD_5\splus2\bA_4&
 \ref{fig.e6-h}(f)&(1,0)
 &\tref{s.group.2loop}&l=5\cr
\exstrut&&&&&\cr
\noalign{\hrule}
\exstrut&&&&&\cr
\2 39&\bE_6\splus\bA_5\splus4\bA_2&
 \ref{fig.hex}(a)&(1,0)
 &\tiref{s.group.hexagon}{39}&(3,3,6,3,3,\0)\cr
\1 40&\same9 (\bE_6\splus\bA_8\splus\bA_2)\splus\bA_2\splus\bA_1\!\!\!&
 \ref{fig.hex}(b)&(0,1)
 &\tiref{s.group.hexagon}{40}&(9,\0,3,3,\0,2)\cr
\2 41&\same3 (\bE_6\splus\bA_5\splus2\bA_2)\splus\bA_4&
 \ref{fig.hex}(c)&(1,0)
 &\tiref{s.group.hexagon}{41}&(6,\0,5,3,3,\0)\cr
42&\same13 \bE_6\splus\bA_6\splus\bA_4\splus\bA_2\splus\bA_1&
 \ref{fig.hex}(d)&(0,1)
 &\tref{s.group.hexagon}&(7,\0,3,5,2,\0)\cr
\exstrut&&&&&\cr
\noalign{\hrule}
\crcr}}}
\endinsert

\midinsert
\table\label{tab.e6-r}
Maximal sets of singularities with a type~$\bE_6$ point
represented by reducible sextics
\endtable
\deftable

\def\no{}
\def\FRAG{fig.e6-r}
\def\tref#1{\,\text{\ref{#1}}\,}
\def\tiref#1#2{\,\text{\iref{#1}{#2}}\,}
\def\same{\afterassignment\dosame\count0=}
\def\dosame{\llap{$^{\the\count0\,}$}}
\centerline{\vbox{\offinterlineskip\halign{%
\tabstrut\ \ \ \hss#\rlap{$'$\hss}\hss\ \vrule&
\quad\,$#$\hss\quad\vrule&
\ \,\null#\hss\ \vrule&&\,\hss$#$\hss\,\vrule\cr
\noalign{\hrule}
\exstrut&&&&&\cr
\omit\tabstrut\ \#\ \vrule&\text{Set of singularities}&
 \hss Figure&\text{Count}&\pi_1&\text{Parameters}\cr
\exstrut&&&&&\cr
\noalign{\hrule}
\exstrut&&&&&\cr
\2 1&\same14 (\bE_6\splus2\bA_5)\splus\bA_3&
 \fragi(a)1&(1,0)&\tiref{s.group.loop-r}{1'}&(6,4,6,\0)\cr
2&\same15 \bE_6\splus\bA_7\splus\bA_5\splus\bA_1&
 \bifragi(a)2&(0,1)&\tref{s.group.loop-r}&(\0,\0,2,\0)\cr
\1 3&\same16 \bE_6\splus\bA_7\splus\bA_3\splus\bA_2\splus\bA_1&
 \fragi(a)3&(1,0)&\tiref{s.group.loop-r}{3'}&(4,8,\0,3)\cr
4&\same17 \bE_6\splus\bA_6\splus\bA_5\splus2\bA_1&
 \fragi(b)1&(1,0)&\tref{s.group.loop-r}&(\0,\0,2,\0)\cr
5&\bE_6\splus\bA_6\splus2\bA_3\splus\bA_1&
 \fragi(b)2&(1,0)&\tref{s.group.loop-r}&(7,4,4,\0)\cr
\1 6&\same18 \bE_6\splus\bA_5\splus\bA_4\splus\bA_3\splus\bA_1&
 \fragi(b)3&(1,0)&\tiref{s.group.loop-r}{6'}&(5,6,4,\0)\cr
7&\same10 \bE_6\splus\bA_9\splus\bA_4&
 \fragi(c)1&(1,0)&\tref{s.group.loop-r}&(\0,\0,\0,1)\cr
8&\same19 \bE_6\splus\bA_9\splus\bA_3\splus\bA_1&
 \fragi(c)2&(1,0)&\tref{s.group.loop-r}&(4,10,\0,2)\cr
9&\bE_6\splus\bD_9\splus\bA_3\splus\bA_1&
 \ref{fig.e6-r}(d)&(1,0)&\tref{s.group.loop-r}&(\0,\0,2,\0)\cr
10&\bE_6\splus\bD_8\splus\bA_4\splus\bA_1&
 \ref{fig.e6-r}(d)&(1,0)&\tref{s.group.loop-r}&(\0,\0,2,\0)\cr
11&\bE_6\splus\bD_6\splus\bA_4\splus\bA_3&
 \ref{fig.e6-r}(d)&(1,0)&\tref{s.group.loop-r}&(5,4,\0,\0)\cr
12&\bE_6\splus\bD_5\splus\bA_4\splus\bA_3\splus\bA_1&
 \ref{fig.e6-r}(d)&(1,0)&\tref{s.loop.D}&(\0,\0,2,\0)\cr
13&\bE_6\splus\bD_{10}\splus\bA_3&
 \ref{fig.e6-r}(e)&(1,0)&\tref{s.group.loop-r}&(\0,\0,\0,1)\cr
14&\bE_6\splus\bD_8\splus\bA_5&
 \ref{fig.e6-r}(e)&(1,0)&\tref{s.group.loop-r}&(\0,\0,\0,1)\cr
\1 15&\same20 \bE_6\splus\bD_5\splus\bA_5\splus\bA_3&
 \ref{fig.e6-r}(e)&(1,0)&\tiref{s.group.loop-r}{15'}&(4,6,\0,\0)\cr
\1 16&\same20 \bE_6\splus\bD_5\splus\bA_5\splus\bA_3&
 \ref{fig.e6-r}(e)&(1,0)&\tref{s.loop.D}&(4,6,\0,1)\cr
17&\bE_6\splus\bD_7\splus\bD_6&
 \ref{fig.e6-r}(f)&(1,0)&\tref{s.loop.noA}&\cr
18&\bE_6\splus\bD_7\splus\bD_5\splus\bA_1&
 \ref{fig.e6-r}(f)&(1,0)&\tref{s.loop.D}&(\0,2,\0,\0)\cr
19&\bE_6\splus\bD_6\splus\bD_5\splus\bA_2&
 \ref{fig.e6-r}(f)&(1,0)&\tref{s.loop.D}&(3,\0,\0,\0)\cr
\exstrut&&&&&\cr
\noalign{\hrule}
\exstrut&&&&&\cr
\2 20&(\bE_6\splus\bA_5\splus2\bA_2)\splus\bA_3\splus\bA_1\!\!\!&
 \ref{fig.hex-r}(a)&(1,0)&\tiref{s.group.hex.1}{20'}&(2,6,3,4,3,\0)\cr
\4 21&\bE_6\splus\bA_4\splus2\bA_3\splus\bA_2\splus\bA_1&
 \ref{fig.hex-r}(b)&(1,0)&\tiref{s.group.hex.1}{21'}&(2,5,4,3,4,\0)\cr
\1 22&\bE_6\splus\bA_7\splus\bA_3\splus3\bA_1&
 \ref{fig.hex-r}(c)&(1,0)&\tiref{s.AM}{22'}&(2,8,2,\0,4,\0)\cr
\1 23&\bE_6\splus\bA_5\splus2\bA_3\splus2\bA_1&
 \ref{fig.hex-r}(d)&(1,0)&\tiref{s.AM}{23'}&(2,4,6,4,\0,\0)\cr
\2 24&(\bE_6\splus2\bA_5)\splus3\bA_1&
 \ref{fig.hex-r}(e)&(1,0)&\tiref{s.AM}{24'}&(2,6,6,2,\0,\0)\cr
25&\same19 \bE_6\splus\bA_9\splus\bA_3\splus\bA_1&
 \ref{fig.hex-r}(f)&(1,0)&\tref{s.group.hex.2}&(10,\0,4,2,\0,\0)\cr
26&\same21 \bE_6\splus\bA_{11}\splus2\bA_1&
 \ref{fig.hex-r}(g)&(1,0)&\tref{s.group.hex.2}&(12,\0,2,\0,2,\0)\cr
27&\same15 \bE_6\splus\bA_7\splus\bA_5\splus\bA_1&
 \ref{fig.hex-r}(h)&(0,1)&\tref{s.group.hex.2}&(8,\0,2,\0,6,\0)\cr
\2 28&\same21 (\bE_6\splus\bA_{11})\splus2\bA_1&
 \ref{fig.hex-r}(i)&(1,0)&\tref{s.group.hex.2}&(12,\0,2,\0,\0,2)\cr
\2 29&\same14 (\bE_6\splus2\bA_5)\splus\bA_3&
 \ref{fig.hex-r}(j)&(1,0)&\tref{s.group.hex.2}&(6,\0,4,6,\0,\0)\cr
\1 30&\same18 \bE_6\splus\bA_5\splus\bA_4\splus\bA_3\splus\bA_1&
 \ref{fig.hex-r}(k)&(1,0)&\tiref{s.group.hex.1}{30'}&(6,5,4,2,\0,\0)\cr
31&\same17 \bE_6\splus\bA_6\splus\bA_5\splus2\bA_1&
 \ref{fig.hex-r}(l)&(1,0)&\tref{s.group.hex.1}&(6,7,2,\0,2,\0)\cr
32&\bE_6\splus\bA_7\splus\bA_4\splus2\bA_1&
 \ref{fig.hex-r}(m)&(1,0)&\tref{s.group.hex.1}&(8,5,2,\0,\0,2)\cr
\1 33&\same16 \bE_6\splus\bA_7\splus\bA_3\splus\bA_2\splus\bA_1&
 \ref{fig.hex-r}(n)&(1,0)&\tiref{s.group.hex.1}{33'}&(8,3,\0,4,2,\0)\cr
\1 34&\bE_6\splus\bA_9\splus\bA_2\splus2\bA_1&
 \ref{fig.hex-r}(o)&(1,0)&\tiref{s.group.hex.1}{34'}&(10,3,\0,2,\0,2)\cr
35&\bE_6\splus\bD_8\splus\bA_3\splus2\bA_1&
 \ref{fig.hex-r}(p)&(1,0)&\tref{s.group.hex.D}&(\0,2,\0,4,2,\0)\cr
\1 36&\bE_6\splus\bD_6\splus2\bA_3\splus\bA_1&
 \ref{fig.hex-r}(p)&(1,0)&\tiref{s.AM}{36'}&(2,4,4,\0,\0,\0)\cr
\1 37&\bE_6\splus\bD_{10}\splus3\bA_1&
 \ref{fig.hex-r}(q)&(1,0)&\tiref{s.AM}{37'}&(2,\0,2,\0,2,\0)\cr
38&\bE_6\splus\bD_6\splus\bA_5\splus2\bA_1&
 \ref{fig.hex-r}(q)&(1,0)&\tref{s.group.hex.D}&(\0,6,2,\0,2,\0)\cr
\exstrut&&&&&\cr
\noalign{\hrule}
\crcr}}}
\endinsert

\Remark\label{rem.4,5}
For pairs \Nos$4$, $5$ and \Nos$6$, $7$, one can ask if the two
curves remain non-equivalent if a permutation of the two
type~$\bE_6$ points is allowed. For \Nos$6$ and~$7$, there still are
two distinct equisingular
deformation families, see~\cite{degt.e6}
or~\cite{OkaPho}; for \Nos$4$ and~$5$, the two curves become
equivalent, see~\cite{Shimada}.
(Alternatively, if \Nos$4$ and~$5$ were not equivalent, each of
the curves
would have a symmetry interchanging its two type~$\bE_6$ points.
Since the curve is maximizing, the symmetry would necessarily be
stable, contradicting to~\cite{symmetric}.)
\endRemark

\Remark
The sets of singularities \Nos$3$ and~$8$ with $(n_r,n_c)=(0,1)$
can be realized by real
curves, see~\cite{Shimada},
but with respect to this real structure the two
type~$\bE_6$ points must be complex conjugate.
\endRemark

\section{The computation\label{S.computation}}

\subsection{Preliminaries and notation}\label{s.group}
To compute the fundamental groups, we use Zariski--van Kampen's
method~\cite{vanKampen}, applying it to the ruling of~$\Sigma_4$.
The
principal steps of the computation are outlined below; for more
details, see~\cite{dessin.e7} and~\cite{dessin.e8}.

Fix a maximizing sextic~$B$ satisfying~$(*)$ and let~$\B$ be the
maximal trigonal curve given by Proposition~\ref{1-1.e6}.
We take for the fiber at infinity~$F_\infty$ the distinguished
type~$\tA{_5}$ fiber~$\F$ of~$\B$, and
for the reference fiber~$F$, the fiber over an appropriate
vertex $v$ of the skeleton~$\Sk$ of~$\B$
in the boundary~$\partial\H$ of
the hexagonal region~$\H$ containing~$\F$. Choose a marking
at~$v$, see~\cite{degt.kplets}, so that
the edges~$e_2$ and~$e_3$ at~$v$
belong to the boundary~$\partial\H$, and let
$\{\Ga_1,\Ga_2,\Ga_3\}$ be a canonical basis in~$F$ defined by
this marking, see~\cite{degt.kplets} or Figure~\ref{fig.basis}.
(The precise choice of the vertex~$v$ and the marking is explained
below on a case by case basis.) Denote $\Gr=\Ga_1\Ga_2\Ga_3$.

\midinsert
\centerline{\picture{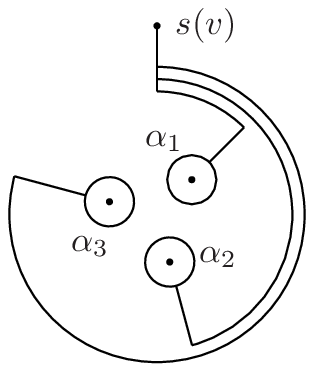}}
\figure
A canonical basis $\{\Ga_1,\Ga_2,\Ga_3\}$
\endfigure\label{fig.basis}
\endinsert

According to~\cite{dessin.e7}, the generators $\Ga_1$, $\Ga_2$,
$\Ga_3$ are subject to the so called \emph{relation at infinity}
$$
\Gr^4=(\Ga_2\Ga_3)^3.
\eqtag\label{rel.inf}
$$
Let $F_1,\ldots,F_r$ be the singular fibers of~$\B$ other
than~$\F$.
Dragging~$F$ about~$F_j$
and keeping the base point in an appropriate section, one obtains
an automorphism $m_j\in\BG3\subset\Aut\<\Ga_1,\Ga_2,\Ga_3\>$,
called the
\emph{braid monodromy} about~$F_j$.
In this notation, the group $\pi_1:=\pi_1(\Cp2\sminus B)$ has a
presentation of the form
$$
\pi_1=\bigl<\Ga_1,\Ga_2,\Ga_3\bigm|
 \text{$m_j=\id$, $j=1,\ldots,r$, and
 \eqref{rel.inf}}\bigr>,
\eqtag\label{eq.vanKampen}
$$
where each \emph{braid relation} $m_j=\id$, $j=1,\ldots,r$,
should be understood as
the triple of relations $m_j(\Ga_i)=\Ga_i$, $i=1,2,3$.
Furthermore, in the presence of the relation at infinity, (any)
one of the braid relations $m_j=\id$, $j=1,\ldots,r$, can be
omitted.

The braid monodromies~$m_j$ are computed using~\cite{degt.kplets};
all necessary details are explained below in this section.

Throughout the paper, all finite
groups/quotients are analyzed using \GAP~\cite{GAP}; most infinite
groups are handled by means of the following obvious lemma, which we
state here for future references.

\lemma\label{center}
Let~$G$ be a group, and let $a\in G$ be a central element whose
projection to the abelianization $G/[G,G]$ has infinite order.
Then the projection $G\to G/a$ restricts to an isomorphism
$[G,G]=[G/a,G/a]$.
\qed
\endlemma

Given two
elements~$\Ga$, $\Gb$ of a group and a nonnegative integer~$m$,
introduce the notation
$$
\{\Ga,\Gb\}_m=\cases
(\Ga\Gb)^k(\Gb\Ga)^{-k},&\text{if $m=2k$ is even},\\
\bigl((\Ga\Gb)^k\Ga\bigr)\bigl((\Gb\Ga)^k\Gb\bigr)\1,&
 \text{if $m=2k+1$ is odd}.
\endcases
$$
The relation $\{\Ga,\Gb\}_m=1$ is equivalent to $\Gs^m=\id$, where
$\Gs$ is the generator of the braid group~$\BG2$ acting on the
free group $\<\Ga,\Gb\>$. Hence,
$$
\{\Ga,\Gb\}_m=\{\Ga,\Gb\}_n=1\quad
\text{is equivalent to}\quad
\{\Ga,\Gb\}_{\gcd(m,n)}=1.
\eqtag\label{eq.equiv}
$$
For the small values of~$m$, the relation $\{\Ga,\Gb\}_m=1$ takes
the following form:
\Dashes
\dash
$m=0$: tautology;
\dash
$m=1$: the identification $\Ga=\Gb$;
\dash
$m=2$: the commutativity relation $[\Ga,\Gb]=1$;
\dash
$m=3$: the braid relation $\Ga\Gb\Ga=\Gb\Ga\Gb$.
\endDashes

\subsection{Two type~$\bE_6$ singular points}\label{s.group.2e6}
It suffices to consider the set of singularities
$2\bE_6\splus\bA_6\splus\bA_1$ (\No$8$ in Table~\ref{tab.e6})
only. The fundamental groups of all sextics of torus type with two
type~$\bE_6$ singular points are found in~\cite{degt.e6}, and the
remaining curves that are not of torus type are covered
by~\ref{s.group.loop} and~\ref{s.group.2loop} below.

Thus, consider the skeleton~$\Sk$
given by Figure~\ref{fig.+e6}(c)--$3$.
Let~$u$ be the monovalent \black-vertex of~$\Sk$,
and let~$v$ be the
trivalent \black-vertex adjacent to~$u$.
Mark~$v$ so that
$[u,v]$ is the edge~$e_2$
at~$v$. Then, in addition to~\eqref{rel.inf}, the group~$\pi_1$
has the
relations
$$
\{\Ga_2\Ga_3\Ga_2\1,\Ga_1\}_7=
\{\Ga_1\Ga_2\Ga_3\Ga_2\1\Ga_1\1,\Ga_2\}_2=1,\qquad
\Gr\Ga_2\Gr\1=\Ga_2\Ga_3\Ga_2\1,
$$
obtained from the heptagon, the bigon, and the monovalent
\black-vertex~$u$ of~$\Sk$, respectively.
The resulting group is abelian.

\subsection{Hexagon with a loop: irreducible curves}\label{s.group.loop}
Assume that $\H$ is a hexagon with a loop, see~\ref{s.hloop} and
Figure~\ref{fig.regions}, and take for~$v$ the vertex shown in the
figure.

\midinsert
\centerline{\cpic{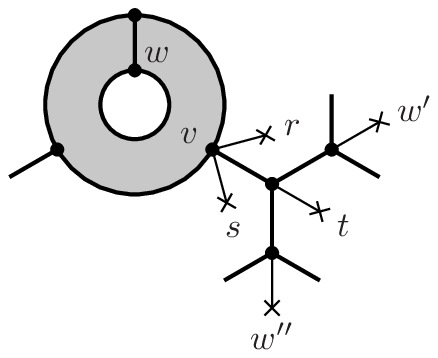}}
\figure
A hexagon with a loop: the regions
\endfigure\label{fig.regions}
\endinsert

The inner loop of~$\H$ (the monogonal region containing~$w$ in
its boundary) gives the relation
$$
\Ga_1=\Ga_3\1\Ga_2\Ga_3.
\eqtag\label{rel.inloop}
$$
Extend~$\Sk$ to a dessin, see~\cite{degt.kplets}, and consider the
\cross-vertices~$r$, $s$, $t$, and~$w'$ shown in
Figure~\ref{fig.regions}. (We do not assert that all these
vertices are distinct.) Assume that they are at the centers of
$l$-, $m$-, $n$-, and $k$-gonal regions of~$\Sk$, respectively.
Then the braid relations about the singular fibers of~$\B$ over
these vertices are
$$
\aligned
r:&\quad\{\Ga_1,\Ga_2\}_l=1,\\
s:&\quad\{\Ga_1,\Ga_2\Ga_3\Ga_2\1\}_m=1,\\
t:&\quad\{\Ga_2,\Gr\Ga_3\Gr\1\}_n=1,\\
w':&\quad\{\Ga_2\1\Ga_1\Ga_2,\Gr\Ga_3\Gr\1\}_k=1,
\endaligned\eqtag\label{rel.loop}
$$
assuming that the fibers are of type~$\tA{}$.
For a $\tD{}$ type fiber,
we omit the corresponding relation
in~\eqref{rel.loop} and indicate this fact by a `$-$' in the
parameter list.
(Sometimes, we also omit a relation just because it is not
necessary to prove that $\pi_1$ is abelian.)
Using the values of $(l,m,n,k)$ shown in
Table~\ref{tab.e6}, one concludes that the groups of most sextics
that are not of torus type are abelian. The few exceptional cases
are treated separately below.

The same arguments apply to the sets of singularities
$2\bE_6\splus\bA_7$ (\No$3$ in Table~\ref{tab.e6}) and
$2\bE_6\splus\bA_4\splus\bA_3$ (\No$5$ in Table~\ref{tab.e6}), as
the corresponding skeletons can be represented as shown in
Figure~\ref{fig.alt}(a)--$1,\bar1$ and
Figure~\ref{fig.alt}(a)--$2$, respectively. The former
fundamental group is
abelian, the latter is the order~$720$ group given by
$$
\pi_1=\<\Ga_1,\Ga_2,\Ga_3\,|\,
 \text{\eqref{rel.inf}, \eqref{rel.inloop}, \eqref{rel.loop}}\>
\eqtag\label{rep.hloop}
$$
with $(l,m,n,k)=(4,5,\0,\0)$; it splits into semidirect product
$\SL(2,\FF_5)\rtimes\CG6$.
(Recall that the braid relation about
the remaining type~$\tE{_6}$ singular fiber can be ignored.)
An alternative presentation of this group is given
by~\eqref{rep.h2loop}.

\midinsert
\centerline{\vbox{\halign{\hss#\hss&&\qquad\qquad\qquad\hss#\hss\cr
\cpic{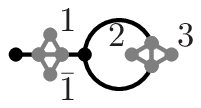}&\cpic{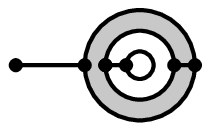}\cr
\noalign{\medskip}
(a)&(b)\cr
\crcr}}}
\figure
\endfigure\label{fig.alt}
\endinsert

\subsubsection{The set of singularities
$\bE_6\splus\bA_{10}\splus\bA_3$ \rom(\No$11$ in
Table~\ref{tab.e6}\rom)}\label{s.e6+a10+a3}
In addition to \eqref{rel.inf}, \eqref{rel.inloop},
and~\eqref{rel.loop} with $(l,m,n,k)=(4,11,\0,\0)$, one also
has the relation
$$
(\Ga_1\Ga_2\Ga_3\Ga_2\1)\Ga_1(\Ga_1\Ga_2\Ga_3\Ga_2\1)\1=
 (\Ga_2\1\Ga_1\Ga_2)\Ga_3(\Ga_2\1\Ga_1\Ga_2)\1
$$
from the lower left loop in Figure~\ref{fig.e6}(b). The resulting
group is abelian.

Alternatively,
choosing a canonical basis $\{\Ga_1,\Ga_2,\Ga_3\}$ in the fiber
over the upper left \black-vertex in Figure~\ref{fig.e6}(b), one
obtains the relations $\Ga_2=\Ga_3$ (from the upper left loop) and
$\Ga_2\1\Ga_1\Ga_2=\Gr\1\Ga_2\Gr$ (from the lower left loop). In
view of the former, the latter simplifies to the braid relation
$\Ga_1\Ga_2\Ga_1=\Ga_2\Ga_1\Ga_2$, \ie, $\{\Ga_1,\Ga_2\}_3=1$. On
the other hand, the $11$-gonal outer region of the skeleton gives
the relation $\{\Ga_1,\Ga_2\}_{11}=1$. Hence
$\{\Ga_1,\Ga_2\}_1=1$, see~\eqref{eq.equiv},
and the group is cyclic.

\subsubsection{The set of singularities
$\bE_6\splus\bD_5\splus\bA_6\splus\bA_2$ \rom(\No$32$ in
Table~\ref{tab.e6}\rom)}\label{s.e6+d5+a6+a2}
In this case we assume that the type~$\tD{_5}$
singular fiber is chosen inside the inner loop of the insertion,
see Remark~\ref{rem.27-28}. Hence the group has no
relation~\eqref{rel.inloop}. However, relations~\eqref{rel.inf}
and~\eqref{rel.loop} with $(l,m,n,k)=(7,3,\0,1)$ suffice to show
that the group is abelian.

\subsubsection{The set of singularities
$\bE_6\splus\bD_5\splus\bA_8$ \rom(\No$28$ in
Table~\ref{tab.e6}\rom)}\label{s.e6+d5+a8}
As above, the type~$\tD{_5}$
fiber is chosen inside the inner loop of the insertion,
see Remark~\ref{rem.27-28}. Hence the group has no
relation~\eqref{rel.inloop}.
Still, it has
relations~\eqref{rel.inf}
and \eqref{rel.loop} with $(l,m,n,k)=(9,\0,1,\0)$ and, in
addition, the relation
$$
\Ga_2\Ga_3\Ga_2\1=(\Ga_1\Ga_2\Ga_3)\Ga_2(\Ga_1\Ga_2\Ga_3)\1
$$
resulting from the left loop in Figure~\ref{fig.e6}(f).
These relations suffice to show that the group is abelian.

\subsubsection{The set of singularities
$(\bE_6\splus\bA_{11})\splus\bA_2$ \rom(\No$12$ in
Table~\ref{tab.e6}\rom)}\label{s.e6+a11+a2}
As explained in~\ref{s.e6+a10+a3},
the group~$\pi_1$ is a quotient
of the braid group~$\BG3$ (for the braid relations, only the two
left loops in Figure~\ref{fig.e6}(b) are used);
hence $\pi_1=\bar\BG3$ due to the following simple lemma (see,
\eg, Lemma~3.6.1 in~\cite{degt.2a8}).

\lemma\label{lem.BG3}
Let $B$ be an irreducible plane sextic of torus type. Then any
epimorphism $\BG3\onto\pi_1(\Cp2\sminus B)$ factors through an
isomorphism $\bar\BG3=\pi_1(\Cp2\sminus B)$.
\endlemma

\subsubsection{The set of singularities
$(\bE_6\splus\bA_5\splus2\bA_2)\splus\bA_4$ \rom(\No$9$ in
Table~\ref{tab.e6}\rom)}\label{s.e6+a5+2a2+a4}
\mnote{Former Subsubsection 3.3.9 is split into two and edited}
The group $\pi_1$ has presentation~\eqref{rep.hloop} with
$(l,m,n,k)=(6,5,3,\0)$. Using \GAP~\cite{GAP}, one can see that
the
quotient $\pi_1/\Ga_1^5$ is a perfect group of order~$7680$,
whereas $\RBG3/\Gs_1^5=\AG5$ has order~$60$.
Hence the natural epimorphism
$\pi_1\onto\RBG3$ is proper.

\subsubsection{Other curves of torus type}\label{s.loop.torus}
The remaining sextics of torus type appearing in this
section are
\Nos$13$ and~$18$ in Table~\ref{tab.e6}.
The group~$\pi_1$ has presentation~\eqref{rep.hloop} with the
values of the parameters $(l,m,n,k)$ given in the table.
Both groups factor to~$\bar\BG3$. Entering the presentations to
\GAP~\cite{GAP} and simplifying them \via
\smallskip
\halign{\qquad\tt#\hss\cr
P := PresentationNormalClosure(g, Subgroup(g, [g.1/g.2]));\cr
SimplifyPresentation(P);\cr}
\smallskip\noindent
one finds that
the commutant $[\pi_1,\pi_1]$ is a
free group on two generators. Since all groups involved are
residually finite, hence Hopfian, one concludes that both
epimorphisms $\pi_1\onto\bar\BG3$ are isomorphisms. (This approach
was suggested to me by E.~Artal Bartolo.)


\subsection{Hexagon with a loop: reducible curves}\label{s.group.loop-r}
Choose a basis $\{\Ga_1,\Ga_2,\Ga_3\}$ as in
Subsection~\ref{s.group.loop}.
Then \eqref{rel.inf} and~\eqref{rel.inloop} imply that
$\pi_1/[\pi_1,\pi_1]=\Z$ and the projection $\pi_1\to\Z$ is
given by $\Ga_1,\Ga_2\mapsto1$, $\Ga_3\mapsto-5$. In particular,
it follows that the sextic splits into an irreducible
quintic and a line.

In addition to~\eqref{rel.inf} and~\eqref{rel.inloop},
consider the relations
$$
\aligned
r:&\quad\{\Ga_1,\Ga_2\}_l=1,\\
s:&\quad\{\Ga_1,\Ga_2\Ga_3\Ga_2\1\}_m=1,\\
t:&\quad\{\Ga_2,\Gr\Ga_3\Gr\1\}_n=1,\\
w'':&\quad\{\Ga_2\1\Ga_1\Ga_2,\Gr\1\Ga_2\Gr\}_k=1
\endaligned\eqtag\label{rel.loop-r}
$$
arising from the \cross-vertices $r$, $s$, $t$, $w''$ in
Figure~\ref{fig.regions} (assuming that the fibers over these
vertices
are of type~$\tA{}$; if a fiber is of type~$\tD{}$, the
corresponding relation is omitted.)

In order to analyze the group using \GAP~\cite{GAP}, observe that
relation~\eqref{rel.inf} implies $[\Ga_1,(\Ga_2\Ga_3)^3]=1$; then,
in view of~\eqref{rel.inloop}, the element $(\Ga_2\Ga_3)^3$
commutes with~$\Ga_3$ and hence with~$\Ga_2$. Thus,
$(\Ga_2\Ga_3)^3\in\pi_1$ is a central element and
its projection
to the
abelianization of~$\pi_1$ is the element $-12$
of infinite order. Due to
Lemma~\ref{center}, it suffices to study the commutant of the
quotient $\pi_1/(\Ga_2\Ga_3)^3$. The abelianization of the latter
quotient is $\CG{12}$.

The sets of parameters $(l,m,n,k)$ used in the calculation are
listed in Table~\ref{tab.e6-r}.
The curves with the following sets of singularities
have nonabelian groups:
\roster
\item\local{1'}
$(\bE_6\splus2\bA_5)\splus\bA_3$ (\No$1'$):
the curve is of torus type, hence
$\ord[\pi_1,\pi_1]=\infty$. Note that $\pi_1/\Ga_1^2=\GL(2,\FF_3)$
has order $48$; in
particular, $\pi_1$ is not isomorphic to any braid group~$\BG{n}$.
\item\local{3'}
$\bE_6\splus\bA_7\splus\bA_3\splus\bA_2\splus\bA_1$ (\No$3'$):
one has $\pi_1=\SL(2,\FF_7)\times\Z$.
\item\local{6'}
$\bE_6\splus\bA_5\splus\bA_4\splus\bA_3\splus\bA_1$ (\No$6'$):
one has $\pi_1=\SL(2,\FF_5)\rtimes\Z$.
\item\local{15'}
$\bE_6\splus\bD_5\splus\bA_5\splus\bA_3$ (\No$15'$):
one has $\pi_1=((\CG3\times\CG3)\rtimes\CG3)\rtimes\Z$.
\endroster
(In Items~\loccit{3'}, \loccit{6'}, and~\loccit{15'},
the centralizer of
$[\pi_1,\pi_1]$ projects to a subgroup of index~$1$, $2$, and~$4$,
respectively, in the
abelianization. In the former case, it follows that the product
is direct.) The group in Item~\loccit{3'} was first computed
in~\cite{Artal.braids}.

\subsubsection{The set of singularities
$\bE_6\splus\bD_7\splus\bD_6$
\rm(\No$17'$ in Table~\ref{tab.e6-r})}\label{s.loop.noA}
The curve has no $\tA{}$-type singular fibers outside the
insertion, and we replace~\eqref{rel.loop-r} with the relations
$$
(\Gr\Ga_1\Ga_2\Ga_1)\Ga_2(\Gr\Ga_1\Ga_2\Ga_1)\1=\Ga_1,\qquad
 (\Gr\Ga_1\Ga_2)\Ga_1(\Gr\Ga_1\Ga_2)\1=\Ga_2
$$
resulting from the type~$\tD{_7}$ singular fiber over~$r$. The
resulting group is abelian.

\subsubsection{A $\tD{}$ type fiber inside the insertion}\label{s.loop.D}
If the singular fiber of~$\B$
in the loop inside the insertion is of type~$\tD{_5}$,
relation~\eqref{rel.inloop} should be replaced with
$$
 \Gr\Gb_2\Gb_3\Gb_2\1\Gr\1=\Gb_2,\qquad
 \Gr\Gb_2\Gr\1=\Gb_3,
\eqtag\label{rel.loop.D}
$$
where $\{\Gb_1,\Gb_2,\Gb_3\}$ is an appropriate canonical basis
over the \black-vertex~$w$ in Figure~\ref{fig.regions}.
Using~\cite{degt.kplets},
one
has $\Gb_1=\Ga_1\Ga_3\Ga_1\1$, $\Gb_2=\Ga_1$, and
$\Gb_3=\Ga_3\1\Ga_2\Ga_3$
(in particular, $\Gb_1\Gb_2\Gb_3=\Gr$).
From \eqref{rel.loop.D} it follows that $\Gd:=\Gr^2\Gb_2\Gb_3$
is a central element of~$\pi_1$;
since the projection of~$\Gd$
to the abelianization of~$\pi_1$ is the element $-4$ of infinite
order, one can use Lemma~\ref{center} and study the commutant of
the quotient $\pi_1/\Gd$.

The sets of parameters $(l,m,n,k)$ used in the calculation are
listed in Table~\ref{tab.e6-r}.
The only nonabelian group in this series is the one corresponding
to the set of singularities
$\bE_6\splus\bD_5\splus\bA_5\splus\bA_3$ (\No$16'$); it can be
represented as a
semidirect product
$((\CG3\times\CG3)\rtimes\CG3)\rtimes\Z$
and is isomorphic to the one described
in~\iref{s.group.loop-r}{15'}.

\subsection{Hexagon with a double loop}\label{s.group.2loop}
Assume that $\H$ is a hexagon with a double loop,
see~\ref{s.h2loop}, and choose for~$v$ the vertex shown in
Figure~\ref{fig.h2loop}, right. We can assume that the singular
fibers inside the insertion are all of type~$\tA{}$, see
Remark~\ref{rem.37,38}.
Then, the braid relations resulting from the inner pentagon and
monogon are
$$
\{\Ga_1,\Ga_3\1\Ga_2\Ga_3\}_5=1,\qquad
\Ga_1\Ga_3\1\Ga_2\Ga_3\Ga_1\1=\Ga_3\1\Ga_2\1\Ga_3\Ga_2\Ga_3.
\eqtag\label{rel.h2loop}
$$
In addition, for all curves except \No$37$ in Table~\ref{tab.e6},
there is a relation
$\{\Ga_1,\Ga_2\}_l=1$, where $l=10$, $7$, $6$, $9$, or~$5$ (in the
order of appearance in Table~\ref{tab.e6}). For \No$37$, one has
the commutativity relation $[\Ga_3,\Ga_1\Ga_2]=1$ from the
$\tD{_9}$ type fiber. In all cases, using \GAP~\cite{GAP}, one
concludes that the group is abelian.

The same arguments apply to the set of singularities
$2\bE_6\splus\bA_4\splus\bA_3$ (\No$4$ in Table~\ref{tab.e6}), as
the corresponding skeleton can be represented as shown in
Figure~\ref{fig.alt}(b) (so that one has $l=4$). The resulting
group has order~$720$, and its presentation is
$$
\pi_1=\<\Ga_1,\Ga_2,\Ga_3\,|\,
 \text{\eqref{rel.inf}, \eqref{rel.h2loop}, $\{\Ga_1,\Ga_2\}_4=1$}\>.
\eqtag\label{rep.h2loop}
$$
This group is isomorphic to~\eqref{rep.hloop}, see
Remark~\ref{rem.4,5}.

\subsection{Genuine hexagon: irreducible curves}\label{s.group.hexagon}
Consider one of the four skeletons shown in Figure~\ref{fig.hex}
and take for~$v$ any vertex in~$\partial\H$.
Let
$v_0=v,v_1,\ldots,v_5$ be the vertices in~$\partial\H$
numbered starting from~$v$ in the \emph{clockwise} direction
(which is the
\emph{counterclockwise} direction in the figures, which represent
the complementary hexagon $S^2\sminus\H$).
Mark each vertex similar to
$v_0=v$ and denote by~$R_i$ the region of~$\Sk$
whose boundary contains the
edges~$e_1$ and~$e_2$ at~$v_i$, $i=0,\ldots,5$. Let~$n_i$ be the
number of vertices in~$\dR_i$, $i=0,\ldots,5$,
\ie, assume that $R_i$ is an $n_i$-gon. Then, in addition to the
common relation at infinity~\eqref{rel.inf},
the group~$\pi_1$ has the relations
$$
\{\Gs_2^{i}(\Ga_2),\Ga_1\}_{n_i}=1,\quad i=0,\ldots,5,
\eqtag\label{rel.hex}
$$
resulting from the singular fibers in~$R_i$. If $R_i$ and~$\H$ are
all but at most one regions of~$\Sk$ (which is always the
case for irreducible curves, see Figure~\ref{fig.hex}), then
\eqref{rel.inf} and~\eqref{rel.hex} form a complete set of
relations for~$\pi_1$.
Furthermore, in the sequence $R_0,\ldots,R_5$, some of the regions
coincide; for each region, it suffices to consider only one
instance in the sequence and ignore the other relations by letting the
corresponding parameters $n_i$ equal~$0$; these relations would
follow from the others.

The values of the parameters $(n_0,\ldots,n_5)$ used in the
calculation are listed in Table~\ref{tab.e6}, and the initial
vertex $v=v_0$ is shown in Figure~\ref{fig.hex} in grey.
For the set of singularities
$\bE_6\splus\bA_6\splus\bA_4\splus\bA_2\splus\bA_1$ (\No$42$ in
Table~\ref{tab.e6}), the resulting group is abelian.
The other three curves are of torus type, hence their groups
factor to $\RBG3$. Below is some information on these groups.
\roster
\item\local{39}
$\bE_6\splus\bA_5\splus4\bA_2$ (\No$39$): the curve is a sextic
of torus type of weight~$8$ in the sense of~\cite{degt.Oka}; hence
$\pi_1$ is much larger than~$\RBG3$: its Alexander module is a
direct sum of \emph{two} copies of $\Z[t]/(t^2-t+1)$. An
alternative presentation of this group is found
in~\cite{degt.8a2}.
\item\local{40}
$(\bE_6\splus\bA_8\splus\bA_2)\splus\bA_2\splus\bA_1$ (\No$40$):
similar to~\ref{s.loop.torus},\mnote{item edited}
one can use \GAP~\cite{GAP} to show that
the natural epimorphism $\pi_1\onto\RBG3$ is an
isomorphism.
\item\local{41}
$(\bE_6\splus\bA_5\splus2\bA_2)\splus\bA_4$ (\No$41$): the
quotient $\pi_1/\Ga_1^5$ is a perfect group of order~$7680$,
whereas $\RBG3/\Gs_1^5=\AG5$. Hence the natural epimorphism
$\pi_1\onto\RBG3$ is proper.
(The values of the parameters actually used
are
$(n_0,\ldots,n_5)=(3,3,5,6,6,5)$.) I do not know whether this
group is isomorphic to the one considered
in~\ref{s.e6+a5+2a2+a4}.
\endroster

\Remark\label{rem.extra}
In Items~\loccit{40} and~\loccit{41}, if the reference fiber is
chosen as shown in Figure~\ref{fig.hex}, the group has also
relation $(\Ga_1\Ga_2)\1\Ga_2(\Ga_1\Ga_2)=\Ga_3$,
see~\eqref{rel.hex.2} below, which may simplify its analysis.
(Formally, this relation follows from the others.)
\endRemark

\subsection{Genuine hexagon:
reducible curves}\label{s.group.hex.1}
The approach of Subsection~\ref{s.group.hexagon}
applies to reducible
curves~$\B$ as well, see Figure~\ref{fig.hex-r},
provided that the skeleton of~$\B$ has at
most one region or $\tD{}$~type singular fiber strictly inside the
complementary hexagon $S^2\sminus\H$, \ie, to all skeletons
except those shown in Figures~\ref{fig.hex-r}(f)--(j). (In the
case of a $\tD{}$-type fiber, Figures~\ref{fig.hex-r}(p) and~(q),
the region~$R_i$ containing this fiber should be ignored
in~\eqref{rel.hex}.)

The first two curves (Figures~\ref{fig.hex-r}(a) and~(b)) do not
seem to have any extra features that would facilitate the study of
their groups. (Each of these curves splits into an irreducible
quintic and a line; hence $\pi_1/[\pi_1,\pi_1]=\Z$.)
The values of the parameters $(n_0,\ldots,n_5)$ are
listed in Table~\ref{tab.e6-r} (assuming that $v=v_0$ is the
vertex shown in the figures in grey), and the resulting groups are
as follows:
\roster
\item\local{20'}
$(\bE_6\splus\bA_5\splus2\bA_2)\splus\bA_3\splus\bA_1$ (\No$20'$):
the curve is of torus type, hence $\pi_1$ factors to the braid
group~$\BG3$; in particular, $[\pi_1,\pi_1]$ is
infinite.
\item\local{21'}
$\bE_6\splus\bA_4\splus2\bA_3\splus\bA_2\splus\bA_1$ (\No$21'$):
the commutant $[\pi_1,\pi_1]$ is perfect (one can compute
the Alexander module $A=0$ similar to Subsection~\ref{s.AM} below);
it appears to be
infinite, but at present I do not even know if it is nontrivial.
\endroster

In the rest of this subsection, we consider the skeletons
with one monogonal region inside $\Cp1\sminus H$, \ie, those shown
in Figures~\ref{fig.hex-r}(k)--(o). Choose for the initial vertex
$v=v_0$ the one shown in the figures in grey. Then, the monogonal
region gives an extra relation
$$
(\Ga_2\Ga_1\Ga_2)\1\Ga_1(\Ga_2\Ga_1\Ga_2)=\Ga_3.
\eqtag\label{rel.hex.1}
$$
(Strictly speaking, this relation follows from the others, but its
presence simplifies the calculations. In particular, since each
sextic~$B$ in question is known to be reducible, it follows that
it splits into an irreducible quintic and a line, the projection
$\pi_1\to\pi_1/[\pi_1,\pi_1]=\Z$ being given by
$\Ga_1,\Ga_3\mapsto1$, $\Ga_2\mapsto-5$.) On the other hand,
relation~\eqref{rel.hex} implies
that $\Gd:=(\Ga_1\Ga_2)^{n_0}$
commutes with~$\Ga_1$ and~$\Ga_2$; hence $\Gd$ is a central
element and, using Lemma~\ref{center}, one can study the commutant
of the quotient $\pi_1/\Gd$.

The parameters $(n_0,\ldots,n_5)$ are listed in
Table~\ref{tab.e6-r}. The following three sets of singularities
result in nonabelian fundamental groups.
\roster
\item[3]\local{30'}
$\bE_6\splus\bA_5\splus\bA_4\splus\bA_3\splus\bA_1$ (\No$30'$):
one has $\pi_1=\SL(2,\FF_5)\rtimes\Z$.
\item\local{33'}
$\bE_6\splus\bA_7\splus\bA_3\splus\bA_2\splus\bA_1$ (\No$33'$):
one has $\pi_1=\SL(2,\FF_7)\times\Z$.
\item\local{34'}
$\bE_6\splus\bA_9\splus\bA_2\splus2\bA_1$ (\No$34'$): the
commutant of~$\pi_1$ is a perfect group; it appears infinite,
but I do not know a proof. The commutants of
$\pi_1/\Ga_1^2$ and $\pi_1/\Ga_1^3$ have orders~$60$ and $51840$,
respectively.
\endroster
(In Items~\loccit{30'} and~\loccit{33'},
the centralizer of
$[\pi_1,\pi_1]$ projects to a subgroup of index~$2$ and~$1$,
respectively, in
$\pi_1/[\pi_1,\pi_1]$.
The groups are isomorphic to those
in~\iref{s.group.loop-r}{6'} and~\ditto{3'}, respectively.)
The group in Item~\loccit{33'} was first computed
in~\cite{Artal.braids}, where it is also shown that sextics
\Nos$3'$ and~$33'$ in Table~\ref{tab.e6-r} are Galois conjugate
over~$\Q(\sqrt2)$.

\subsection{Genuine hexagon: two monogons inside $S^2\sminus\H$}\label{s.group.hex.2}
In this subsection, we
consider a skeleton with two monogonal regions strictly inside
$S^2\sminus\H$, \ie, one of those shown in
Figures~\ref{fig.hex-r}(f)--(j). Take for $v=v_0$ the vertex shown
in the figures in grey. Then, in addition to~\eqref{rel.inf}
and~\eqref{rel.hex}, the group has an extra relation
$$
(\Ga_1\Ga_2)\1\Ga_2(\Ga_1\Ga_2)=\Ga_3
\eqtag\label{rel.hex.2}
$$
resulting from the monogon closest to~$v$. In particular, it
follows that the curve splits into an irreducible quartic and
irreducible conic. Furthermore, since $\Gd:=(\Ga_1\Ga_2)^{n_0}$
commutes with~$\Ga_1$ and~$\Ga_2$, it is a central element and one
can use Lemma~\ref{center} and study the commutant of the quotient
$\pi_1/\Gd$.

The parameters $(n_0,\ldots,n_5)$ are listed in
Table~\ref{tab.e6-r}. For the first three curves, the groups are
abelian. (As a consequence,
the curve defined by
the skeleton in Figure~\ref{fig.hex-r}(g), \No$26'$ in
Table~\ref{tab.e6-r}, is not of torus type.) The other two curves
are of torus type; an alternative way to construct these curves
and to compute their fundamental groups is found
in~\cite{degt.2a8}. (To prove that these curves are of torus type,
one can argue that the existence of such curves is shown
in~\cite{degt.2a8}, and \Nos$28$, $29$ are the only candidates
left.)

\subsection{Genuine hexagon with a $\tD{}$-type fiber}\label{s.group.hex.D}
Consider one of the two skeletons shown in Figure~\ref{fig.hex-r}(p)
or~(q) and choose the initial vertex $v=v_0$ so that $R_0$ is the
region containing the $\tD{}$-type fiber~$F$ of~$\B$.
Then, as above, the defining relations for~$\pi_1$
are~\eqref{rel.inf} and~\eqref{rel.hex}, with the contribution
of~$R_0$ ignored in the latter.

However, we do make use of the region~$R_0$ in order to find
central elements in~$\pi_1$.
Let $n=n_0$.
Then $F$ is of type~$\tD{_{n+4}}$, and the braid relations about~$F$
are
$$
\Ga_3\1\Ga_i\Ga_3=\Gs_1^{n+2}(\Ga_i),\qquad i=1,2.
$$
As a consequence,
$$
[\Ga_3,\Ga_1\Ga_2]=1,\qquad
 \text{hence}\quad
 [\Ga_3,\Gr]=[\Ga_1\Ga_2,\Gr]=1,
\eqtag\label{rel.D}
$$
and
$\Gd:=\Ga_3(\Ga_1\Ga_2)^{1+n/2}$ is a central element of~$\pi_1$.
(Note that in all cases $n$ is even.)
Since
$[\Ga_2\Ga_3,\Gr^4]=1$,
see~\eqref{rel.inf},
from~\eqref{rel.D} one has $[\Ga_2,\Gr^4]=1$
and then
$[\Ga_1,\Gr^4]=1$. Thus, $\Gr^4=(\Ga_2\Ga_3)^3$ is also a central
element of~$\pi_1$, and Lemma~\ref{center} applied twice implies
that the commutant $[\pi_1,\pi_1]$
of~$\pi_1$ is isomorphic to that of the
quotient
$\pi_1/\<\Gd,\Gr^4\>$. (It is worth mentioning that $n\ge2$ and
hence the images of $\Gd$ and~$\Gr$ in the abelianization
of~$\pi_1$ are linearly independent.)

The values of the parameters $(n_0,\ldots,n_5)$ are listed in
Table~\ref{tab.e6-r}. (Recall that the initial vertex~$v_0$ is
determined by the position of the $\tD{}$-type fiber, which
depends on the curve.)
Using \GAP~\cite{GAP}, one concludes that, for the sets of
singularities
$$
\bE_6\splus\bD_8\splus\bA_3\splus2\bA_1
\quad\text{and}\quad
\bE_6\splus\bD_6\splus\bA_5\splus2\bA_1
$$
(\Nos$35'$ and~$38'$ in Table~\ref{tab.e6-r}),
the group~$\pi_1$ is abelian, whereas
for the other two curves (\Nos$36'$ and~$37'$) it has infinite
commutant. For a more precise statement, see Subsection~\ref{s.AM}
below, where we compute the Alexander modules of these and a few
other groups.

\Remark\label{4+1}
The fact that the fundamental groups of the reducible
sextics with the sets
of singularities $\bE_6\splus\bD_6\splus2\bA_3\splus\bA_1$ and
$\bE_6\splus\bD_{10}\splus3\bA_1$ have infinite commutants
can also be explained as
follows. Each curve splits into an irreducible quartic~$B_4$ with
a type~$\bE_6$ singular point and a pair of lines. One of the
lines~$B_1$ either is double tangent to~$B_4$ or has a single
point of $4$-fold intersection with~$B_4$. Hence, even patching
back in the other line (which corresponds to letting one of the
canonical generators of~$\pi_1$ equal to~$1$), one obtains a curve
with large fundamental group (which is, respectively, $\BG3$ or
$\TG{3,4}=\<\Ga,\Gb\,|\,\Ga^3=\Gb^4\>$, see~\cite{groups}).
\endRemark

\subsection{Other sextics with two linear components}\label{s.AM}
With the exception of the two curves mentioned in the previous
section, all maximizing
sextics splitting into an irreducible quartic (with a
type~$\bE_6$ singular point) and two lines have fundamental groups
with infinite commutants. In order
to prove this statement and make it more precise, we compute the
so called
Alexander modules of the groups, see~\cite{Libgober2}.

\definition
Let~$G$ be a group, and let $G'=[G,G]$ be its commutant. The
\emph{Alexander module} of~$G$ is the abelian group $G'\!/[G',G']$
regarded as a $\Z[G/G']$-module \via\ the conjugation action
$(a,x)\mapsto a\1xa$, $a\in G/G'$, $x\in G'\!/[G',G']$.
\enddefinition

Abbreviate $\pi_1'=[\pi_1,\pi_1]$ and denote the Alexander module
of~$\pi_1$ by~$A$.

The sextics in question are represented by the skeletons shown in
Figures~\ref{fig.hex-r}(c)--(e), (p), and~(q), and the defining
relations for~$\pi_1$ are~\eqref{rel.inf} and~\eqref{rel.hex}. (As
usual,
if a $\tD{}$ type fiber is present, the corresponding relation
in~\eqref{rel.hex} should be ignored.) The abelianization
$\pi_1/\pi_1'=\Z\oplus\Z$ is generated by the images~$s$, $t$
of~$\Ga_1$, $\Ga_2$, respectively, and the group ring
$\Z[\pi_1/\pi_1']$ can be identified with the ring
$\LL:=\Z[s,s\1,t,t\1]$ of Laurent polynomials in~$s$, $t$.

Each skeleton in question has a bigonal region not containing
a $\tD{}$-type fiber, and we choose the
initial vertex~$v_0$ so that this region is~$R_0$ (see
the grey vertex in the figures; note that,
for Figures~\ref{fig.hex-r}(p) and~(q),
this choice of~$v_0$ differs from that used in
Subsection~\ref{s.group.hex.D}). Then $[\Ga_1,\Ga_2]=1$ and,
using the Reidemeister--Schreier method (see, \eg,~\cite{MKS}),
one can see that\mnote{more precise reference}
$A$ is generated over~$\LL$ by a single element
$a:=\Ga_1^4\Ga_2\Ga_3$. The relation at infinity~\eqref{rel.inf}
transforms into $(s^{-4}+s^{-8})a=(s^{-3}+s^{-6}+s^{-9})a$, or
$$
Q_\infty(s)a=0,\quad\text{where}\quad
 Q_\infty(s):=(s^2-s+1)(s^4-s^2+1).
\eqtag\label{AM.inf}
$$
Alternatively, \eqref{AM.inf} can be rewritten
in the form
$(s-1)s(s^4+s+1)a=-a$, which means that the multiplication by
$(s-1)$ is invertible in~$A$. For this reason, we cancel the
factor $(s-1)$ in all other relations.

For an integer $m\ge0$, denote $P_m(x)=(x^m-1)/(x-1)$.
(In particular, $P_0\equiv0$ and $P_1\equiv1$.)
Observe that, for curves with two linear components, all
integers~$n_i$ in~\eqref{rel.hex} are even, see
Table~\ref{tab.e6-r}.
A relation
$\{\Gs_2^i(\Ga_2),\Ga_1\}_{2r}=1$ results in the following
relation for~$A$:
$$
\cases
(t-1)P_r(st)P_j(s^4)a=0,
 \quad&\text{if $i=2j$ is even, or}\\
P_r(s^3t)\bigl[(1-s^4t)P_j(s^4)+s^{4j}t\bigr]a=0,
\quad&\text{if $i=2j-1$ is odd}.
\endcases
$$
(Recall that we cancel all factors $(s-1)$.) Using the values of
$(n_0,\ldots,n_5)$ listed in Table~\ref{tab.e6-r},
one arrives at the following Alexander modules.
\roster
\item\local{22'}
$\bE_6\splus\bA_7\splus\bA_3\splus3\bA_1$ (\No$22'$):
$A=\Z[s]/Q_\infty(s)$ and $ta=a$.
\item\local{23'}
$\bE_6\splus\bA_5\splus2\bA_3\splus2\bA_1$ (\No$23'$):
$A=\Z[s]/Q_\infty(s)$ and $ta=-s^{-3}a$.
\item\local{24'}
$(\bE_6\splus2\bA_5)\splus3\bA_1$ (\No$24'$):
$A=\Z[s]/(s^2-s+1)$ and $ta=(1+s^{-4})a$.
\item\local{36'}
$\bE_6\splus\bD_6\splus2\bA_3\splus\bA_1$ (\No$36'$):
$A=\Z[s]/(s^2-s+1)$ and
$ta=-s^{-3}a$.
\item\local{37'}
$\bE_6\splus\bD_{10}\splus3\bA_1$ (\No$37'$):
$A=\Z[s]/Q_\infty(s)$ and $ta=a$.
\endroster
(Note that $Q_\infty\mathbin|(s^{12}-1)$; hence $s$ is invertible in
$\Z[s]/Q_\infty$ and one does not need to consider Laurent
polynomials explicitly.)

\Remark
In Items~\itemref{s.AM}{22'} and~\ditto{37'} above, one has
$n_2=2$, hence $(t-1)a=0$. In Items~\ditto{23'}
and~\ditto{36'}, one has $n_1=4$, hence $(1+s^3t)a=0$. In each
case, $t$ is a Laurent polynomial in~$s$ and one can represent~$A$
as a quotient of the $\LL$-module $\Z[s]/Q_\infty$ (with an
appropriate action of~$t$); then, in most cases, all extra
relations follow from relation at infinity~\eqref{AM.inf}.

Denote $Q(s)=s^2-s+1$ and $R(s)=s^4-s^2+1$, so that
$Q_\infty(s)=Q(s)R(s)$.
In
Item~\ditto{36'}, in addition to~\eqref{AM.inf}, one has the
relation $Q(s)S_1(s)a=0$, where $S_1(s):=(s+1)^2$. Since
$R(s)-(s-1)^2S_1(s)=s^2$ is an invertible element, the relation
ideal (in $\Z[s]$) is generated by a single element $s^2-s+1$.

In
Item~\ditto{24'}, the additional relations are
$$
P_3(s^3t)a=(t-1)P_3(st)a=(1-s^4t+s^4)a=0.
$$
From the last relation, one obtains $ta=(1+s^{-4})a$. Hence, again
$t$ is a polynomial in~$s$ and, substituting $t=1+s^{-4}$
to the first two relations, one has
$$
(s^2-s+1)S_1(s)a=(s^2-s+1)S_2(s)a=0,
$$
where
$$
S_1(s):=s^6+s^5+2s^2+2s+1,\qquad
S_2(s):=s^6+2s^5+2s^4+s+1.
$$
One can easily check that
$$
\multline
s(s^3-s-2)(s^4+s^3+s^2-s+1)R(s)
 -s^2(s^4+s^3+s^2-s+1)S_1(s)\\
 +(s^5+s^4+s^3-s^2+s+1)S_2(s)=1.
\endmultline
$$
Hence again the relation ideal
is generated by a single element $s^2-s+1$.
\endRemark

\Remark
In Items~\itemref{s.AM}{22'} and~\ditto{23'} above, the fact that
the groups are infinite can be explained similar to
Remark~\ref{4+1}: each curve splits into an irreducible
quartic~$B_4$ and a pair of lines,
one of them either being double tangent to~$B_4$ or
having a single
point of $4$-fold intersection with~$B_4$.
\endRemark

\section{Perturbations\label{S.pert}}

We fix a maximizing plane sextic~$B$ satisfying condition~$(*)$
in the introduction
and consider a perturbation~$B'$ of~$B$. Throughout this section,
$\B$ stands for the maximal trigonal curve corresponding to~$B$ \via\
Proposition~\ref{1-1.e6}.

\subsection{Perturbations of the type~$\bE_6$ point~$P$}\label{s.pert.P}
Let $U$ be a Milnor ball about the distinguished type~$\bE_6$
singular point~$P$ of~$B$. The group $\pi_1(\MB\sminus B)$ is
generated by three elements $\Gb_1$, $\Gb_2$, $\Gb_3$ subject to
the relations
$$
\bc_3=(\bc_1\bc_2\bc_3)\bc_2(\bc_1\bc_2\bc_3)\1,\qquad
 \bc_2=(\bc_1\bc_2\bc_3)\bc_1(\bc_1\bc_2\bc_3)\1.
\eqtag\label{rel.e6}
$$
According to~\cite{dessin.e7},
the inclusion homomorphism
$\pi_1(\MB\sminus B)\to\pi_1(\Cp2\sminus B)$ is given by
$$
\bc_1\mapsto\Gr\Ga_1\Gr\1,\qquad
\bc_2\mapsto\Ga_1,\qquad
\bc_3\mapsto\Gr\1\Ga_1\Gr,
$$
where $\{\Ga_1,\Ga_2,\Ga_3\}$ is any basis for
$\pi_1(\Cp2\sminus B)$ as in Subsection~\ref{s.group}.
Note that
a type~$\bE_6$ singularity has an order~$3$ automorphism inducing
the automorphism
$$
\bc_1\mapsto(\bc_1\bc_2)\bc_3(\bc_1\bc_2)\1,\qquad
\bc_2\mapsto\bc_1,\qquad
\bc_3\mapsto\bc_2
$$
of the local fundamental
group. However, on the image in $\pi_1(\Cp2\sminus B)$, this latter
automorphism is inner, \viz. it is induced by the conjugation
$g\mapsto\Gr g\Gr\1$. (For proof, one needs to use the
commutativity relation $[\Ga_1,(\Ga_2\Ga_3)^3]=1$, which follows
from~\eqref{rel.inf}.) Hence, the local automorphisms at~$P$ can
be ignored when studying the extra relations in
$\pi_1(\Cp2\sminus B')$ resulting from a perturbation $B\to B'$.

The fundamental groups $\pi_1(\MB\sminus B')$ of small
perturbations $B\to B'$ are found in~\cite{degt.e6}. They are
as follows:
\roster
\item\local{2a2+a1}
$\bE_6\to2\bA_2\splus\bA_1$: one has
$\pi_1(\MB\sminus B)=\BG4$, the additional relations being
$[\bc_1,\bc_3]=\{\bc_1,\bc_2\}_3=\{\bc_2,\bc_3\}_3=1$;
\item\local{2a2}
$\bE_6\to2\bA_2$ (as a further perturbation of~\loccit{2a2+a1})
and $\bE_6\to\bA_5$:
one has
$\pi_1(\MB\sminus B)=\BG3$, the additional relations being
$\bc_1=\bc_3$ and $\{\bc_1,\bc_2\}_3=1$;
\item\local{others}
all others: one has $\pi_1(\MB\sminus B)=\Z$ and the relations are
$\bc_1=\bc_2=\bc_3$.
\endroster
Note that, in the presence of~\eqref{rel.e6}, the additional
relations in Items~\loccit{2a2+a1} and~\loccit{2a2} follow from
the first relation, \viz. $[\bc_1,\bc_3]=1$ and $\bc_1=\bc_3$,
respectively. The first two perturbations are of \emph{torus
type}, \ie, they preserve torus structures of~$B$; the other
perturbations destroy any torus structure with respect to which
$P$ is an inner singularity.

Combining~\loccit{2a2+a1}--\loccit{others} above with the
inclusion homomorphism, we obtain the following extra relations
for the perturbed group $\pi_1(\Cp2\sminus B')$:
\roster
\item
$\bE_6\to2\bA_2\splus\bA_1$: the extra relation is
$[\Ga_1,\Gr^{-2}\Ga_1\Gr^2]=1$;
\item
$\bE_6\to2\bA_2$ and $\bE_6\to\bA_5$:
the extra relation is $[\Ga_1,\Gr^2]=1$;
\item
all others: the extra relation is $[\Ga_1,\Gr]=1$.
\endroster
As a consequence,
if $P$ is perturbed as in~\loccit{others}, one has
$[\Ga_1,\Ga_2\Ga_3]=1$ and \eqref{rel.inf} becomes
$$
\Ga_1^4\Ga_2\Ga_3=1.
\eqtag\label{rel.inf.p3}
$$
In particular, $\pi_1$ is generated by~$\Ga_1$ and~$\Ga_2$
(or by~$\Ga_1$ and~$\Ga_3$) in this case.

\lemma\label{lem.bigon}
Assume that the distinguished hexagon~$\H$ of~$\Sk$
is genuine, one of the regions~$R_i$, see
Subsection~\ref{s.group.hexagon}, is a bigon, and the point~$P$
is perturbed as in~\itemref{s.pert.P}{others} above. Then the
group $\pi_1(\Cp2\sminus B')$ is abelian.
\endlemma

\proof
Taking~$R_i$ for~$R_0$,
one concludes that $[\Ga_1,\Ga_2]=1$; together
with~\eqref{rel.inf.p3}, this relation implies that the group is
abelian.
\endproof

\Remark
Note that almost all skeletons with a genuine hexagon
have a bigonal region;
exceptions are Figure~\ref{fig.hex}(a), (c) and
Figure~\ref{fig.hex-r}(j).
\endRemark

\subsection{Perturbations of $\bA$ type points}\label{s.A}
Extend the skeleton~$\Sk$ of~$B$ to a dessin by inserting a
\white-vertex at the middle of each edge,
inserting a \cross-vertex~$v_R$ at
the center of each region~$R$, and connecting~$v_R$ to the
vertices in the boundary $\dR$ by appropriate edges,
see~\cite{degt.kplets}. Let~$Q$ be a
type~$\bA_p$ singular point of~$B$; it is located in a
type~$\tA{_p}$ singular fiber~$F$ of~$\B$ at the center of a
certain $(p+1)$-gonal region~$R$ of~$\Sk$.
If $Q$ is perturbed,
$F$ splits into singular fibers of types $\tA{_{s_i}}$
($\tA{_0^*}$ if $s_i=0$), $i=1,\ldots,r$, with $\sum(s_i+1)=p+1$.
Geometrically, the $2(p+1)$-valent \cross-vertex~$v_R$
splits into several \cross-vertices of valencies $2(s_i+1)$.

Assume that the braid
relation in $\pi_1(\Cp2\sminus B)$ resulting
from the region~$R$ above is of the form $\{\Gd_1,\Gd_2\}_{p+1}=1$,
where $\Gd_1$, $\Gd_2$ are certain words in $\Ga_1$, $\Ga_2$,
$\Ga_3$, \cf. \eqref{rel.loop}, \eqref{rel.loop-r},
\eqref{rel.hex}. (In other words, $\Gd_1$, $\Gd_2$ are the first
two elements of a canonical basis over a vertex $u\in\dR$
defined by a marking at~$u$, see~\cite{degt.kplets},
with respect to which $e_1$ and~$e_2$
belong to~$\dR$.)

\lemma\label{lem.A}
In the notation above, the relation for the new group
$\pi_1(\Cp2\sminus B')$ resulting from~$R$ is
$\{\Gd_1,\Gd_2\}_{s}=1$, where $s=\gcd(s_i+1)$, $i=1,\ldots,r$.
\endlemma

\proof
The statement follows from the description of the braid monodromy
found in~\cite{degt.kplets} and from~\eqref{eq.equiv}.
\endproof

\corollary\label{cor.A}
If $Q$ is perturbed, the relation $\{\Gd_1,\Gd_2\}_{p+1}=1$
changes to the relation $\{\Gd_1,\Gd_2\}_{s}=1$, where $s<p+1$ is
a divisor of~$(p+1)$. In particular, one has $s=1$ if $(p+1)$ is a
prime.
\qed
\endcorollary

\corollary\label{cor.A.irr}
If $Q$ is a point of intersection of two components of~$B$
\rom(hence $p$ is odd\rom)
and the
perturbation is to be irreducible in a Milnor ball about~$Q$, then
$s$ in Corollary~\ref{cor.A} is an odd divisor of $(p+1)$.
\qed
\endcorollary

\corollary\label{cor.A.nontorus}
If $B$ is of torus type, $Q$ is an inner singularity
\rom(hence $p+1=0\bmod3$\rom), and the torus structure is to be
destroyed,
then
$s$ in Corollary~\ref{cor.A} is not divisible by~$3$.
\qed
\endcorollary

\subsection{Perturbations of $\bD$ type points}\label{s.D}
Let~$Q$ be a type $\bD_p$, $p\ge5$, singular point of~$B$; it is
located in a type~$\tD{_p}$ singular fiber~$F$ of~$\B$ at the
center of a $(p-4)$-gonal region~$R$ of~$\Sk$.

\lemma\label{lem.D}
If a perturbation $B\to B'$ is irreducible in a Milnor
ball~$\MB_Q$ about~$Q$, the group $\pi_1(\Cp2\sminus B')$ is
abelian.
\endlemma

\proof
Under the assumptions, the group $\pi_1(\MB_Q\sminus B')$ is
abelian, see~\cite{dessin.e8}; on the other hand, the inclusion
homomorphism
$\pi_1(\MB_Q\sminus B')\to\pi_1(\Cp2\sminus B')$ is onto, as $Q$
is a triple point of a trigonal curve.
\endproof

Let~$\Gd_1$, $\Gd_2$ be the first
two elements of a canonical basis over a vertex $u\in\dR$
defined by a marking at~$u$, see~\cite{degt.kplets},
with respect to which $e_1$ and~$e_2$
belong to~$\dR$.
(In other words, if $R$ contained an~$\tA{}$ type fiber,
the resulting braid relation would be
$\{\Gd_1,\Gd_2\}_{p-4}=1$,
\cf. \eqref{rel.loop}, \eqref{rel.loop-r},
\eqref{rel.hex}.)

\lemma\label{lem.D.r}
Assume that a proper perturbation $B\to B'$ is still reducible in
a Milnor ball~$\MB_Q$ about~$Q$. Then, in the notation above, the
new group $\pi_1(\Cp2\sminus B)$ has an extra relation
$\{\Gd_1,\Gd_2\}_{s}=1$ for some integer~$s$, $1\le s\le p-2$.
If $B\cap\MB_Q$ has three components \rom($p$ is even\rom)
but $B'\cap\MB_Q$ has two components, then $s$ is odd.
\endlemma

\proof
The statement follows from the computation of the fundamental
group $\pi_1(\MB_Q\sminus B')$ found in~\cite{dessin.e8}.
\endproof

\subsection{Proof of Theorem~\ref{th.pert.nontorus}}\label{pf.pert.nontorus}
We skip sextics of torus type with two or more type~$\bE_6$
singular points (\Nos$1$, $2$, $6$, and~$7$ in Table~\ref{tab.e6};
they are considered in details in~\cite{degt.e6}) and the
maximizing sextic of weight~$8$ (\No$39$ in Table~\ref{tab.e6};
considered in~\cite{degt.8a2}).

All other perturbations with nonabelian group
$\pi_1:=\pi_1(\Cp2\sminus B)$ can be handled on a case by case basis,
using \GAP~\cite{GAP} and modifying the presentation for~$\pi_1$
found in Section~\ref{S.computation}
according to
Subsection~\ref{s.pert.P},
Corollaries~\ref{cor.A}--\ref{cor.A.nontorus}, and
Lemmas~\ref{lem.D} and~\ref{lem.D.r}. (Recall that, due to
Proposition~5.1.1 in~\cite{degt.8a2}, all singular points of~$B$
can be perturbed arbitrarily and independently.) A few details are
given in~\ref{s.first}--\ref{s.last}
below.
In some cases,
the same approach shows that any proper
perturbation, including reducible ones, that is not of torus type
has abelian fundamental group:
\Dashes

\dash
\No$1'$: all perturbations not of torus type are abelian;

\dash
\No$3'$: all proper perturbations are abelian, see~\ref{s.3'};

\dash
\Nos$6'$, $15'$, $16'$: all proper perturbations are abelian;

\dash
\No$20'$: all perturbations not of torus type are abelian,
see~\ref{s.20'};

\dash
\No$24'$: all perturbations not of torus type are abelian,
see~\ref{s.24'};

\dash
\Nos$28'$, $29'$: all perturbations not of torus type are abelian;

\dash
\Nos$30'$, $33'$: all proper perturbations are abelian;

\dash
\No$34'$: all proper perturbations except
$\bE_6\to2\bA_2\splus\bA_1$ are abelian; for the
(reducible) nonabelian one,
the perturbation epimorphism appears to be an isomorphism.

\endDashes

\subsubsection{The set of singularities
$2\bE_6\splus\bA_4\splus\bA_3$ \rom(\Nos$4$ and~$5$\rom)}\label{s.first}
The two curves are equivalent, see Remark~\ref{rem.4,5}, and it is
easier to use \No$5$ and its model given in
Figure~\ref{fig.alt}(a)--$2$, see Subsection~\ref{s.group.loop}.

\subsubsection{The sets of singularities
$(\bE_6\splus\bA_{11})\splus\bA_2$,
$(\bE_6\splus\bA_8\splus\bA_2)\splus\bA_3$, and
$(\bE_6\splus\bA_8\splus\bA_2)\splus\bA_2\splus\bA_1$
\rom(\Nos$12$, $13$, $18$, and~$40$\rom)}\label{s.12}
One has\mnote{title extended}
$\pi_1=\bar\BG3$, and the statement on
the perturbations follows from Lemma~3.12\mnote{reference updated}
in~\cite{degt.2a8}.

\subsubsection{The set of singularities
$\bE_6\splus\bA_7\splus\bA_3\splus\bA_2\splus\bA_1$ \rom(\No$3'$\rom)}\label{s.3'}
The $\bA_1$ type point is not involved in the computation of the
original group~$\pi_1$, see Subsection~\ref{s.group.loop-r}.
Using Corollary~\ref{cor.A} and computing the braid monodromy as
explained in~\cite{degt.kplets}, one can easily see that the additional
relation resulting from perturbing this point is
$(\Ga_2\1\Ga_1\Ga_2)\Ga_3(\Ga_2\1\Ga_1\Ga_2)\1=\Gr\1\Ga_2\Gr$.

\subsubsection{The set of singularities
$(\bE_6\splus\bA_5\splus2\bA_2)\splus\bA_3\splus\bA_1$
\rom(\No$20'$\rom)}\label{s.20'}
The perturbations of the
distinguished
type~$\bE_6$ point~$P$ are covered by Lemma~\ref{lem.bigon}. If a
type~$\bA_5$ point~$Q$ is perturbed, one can number the regions~$R_i$,
see Subsection~\ref{s.group.hexagon}, so that $R_0$ is a bigon and
$Q$ is over~$R_1$. According to
Corollary~\ref{cor.A.nontorus}, the new group~$\pi_1'$ has relations
$[\Ga_1,\Ga_2]=[\Ga_1,\Ga_3]=1$; then \eqref{rel.inf} turns
into~\eqref{rel.inf.p3} and $\pi_1'$ is generated by two commuting
elements~$\Ga_1$, $\Ga_2$.

Suppose that a type~$\bA_2$ point~$Q$ is perturbed; since
the skeleton is symmetric, one can assume that $Q$ is over
the left triangle in Figure~\ref{fig.hex-r}(a).
Taking this triangle for~$R_0$,
see Subsection~\ref{s.group.hexagon}, one obtains the relations
$\Ga_1=\Ga_2$ (Corollary~\ref{cor.A})
and
$[\Ga_2,\Ga_3\Ga_2\Ga_3\1]=1$ (from the bigon~$R_2$); hence also
$[\Ga_2,\Ga_3\1\Ga_2\Ga_3]=1$. On the other hand,
\eqref{rel.inf} implies $[\Ga_2,(\Ga_2\Ga_3)^3]=1$ (in view of
$\Ga_1=\Ga_2$). Combining these relations, one arrives at
$[\Ga_2,\Ga_3^3]=1$, \ie, $\Ga_3^3$ is a central element.
Now, using
Lemma~\ref{center} and \GAP~\cite{GAP}, one concludes that the
perturbed group is abelian.

\subsubsection{The set of singularities
$(\bE_6\splus2\bA_5)\splus3\bA_1$
\rom(\No$24'$\rom)}\label{s.24'}\label{s.last}
Perturbation of the type~$\bE_6$ or type~$\bA_5$ points are
handled similar to~\ref{s.20'}.
\qed

\subsection{Proof of Theorem~\ref{th.pert.torus}}\label{pf.pert.torus}
As\mnote{title changed}
in the previous subsection, we ignore the sets of
singularities with two or more type~$\bE_6$ points,
see~\cite{degt.e6}, or the set of singularities
$\bE_6\splus\bA_5\splus4\bA_2$ of weight~$8$, see~\cite{degt.8a2}.
These sets of singularities give rise to the exceptions mentioned
in the statement.

Consider a dessin with one of the two fragments shown in
Figure~\ref{fig.fragments}. We assume that a region of the
skeleton may contain several \cross-vertices, of valencies
$2(s_i+1)$, $i=1,\ldots,r$, and parameter~$s$ shown in the figures
stands for the greatest common divisor
$\gcd(s_i+1)$, \cf. Lemma~\ref{lem.A}.

\midinsert
\centerline{\picture{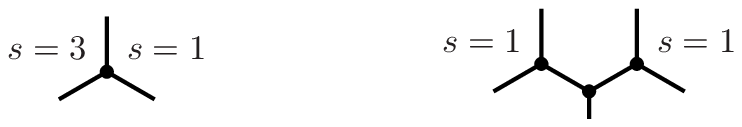}}
\figure
Two special fragments
\endfigure\label{fig.fragments}
\endinsert

\lemma\label{fragments}
Let~$B'$ be an irreducible sextic of torus type and with a
type~$\bE_6$ singular point, and assume that the dessin of the
corresponding trigonal curve~$\B'$ has one of the two fragments
shown in Figure~\ref{fig.fragments}. Then
$\pi_1(\Cp2\sminus B')=\bar\BG3$.
\endlemma

\proof
The statement follows immediately from Lemma~\ref{lem.A}. In the
first case, in an appropriate canonical basis
$\{\Ga_1,\Ga_2,\Ga_3\}$ over the
\black-vertex shown in the figure, one has $\Ga_1=\Ga_2$ and
$\{\Ga_2,\Ga_3\}_3=1$. The second case (the right fragment in the
figure)
is essentially considered
in~\ref{s.e6+a10+a3}. In both cases, one obtains an epimorphism
$\BG3\onto\pi_1(\Cp2\sminus B')$, and Lemma~\ref{lem.BG3} applies.
\endproof

Now, considering sextics of torus type one by one and perturbing
their $\bA$ type points taking into account Corollaries~\ref{cor.A}
and~\ref{cor.A.irr}, one finds only three perturbations
not covered by
Lemma~\ref{fragments}:
\Dashes
\dash
$(\bE_6\splus\bA_5\splus2\bA_2)\splus\bA_4\to
 (\bE_6\splus4\bA_2)\splus\bA_4$
(\No$9$);
\dash
$(\bE_6\splus\bA_8\splus\bA_2)\splus\bA_2\splus\bA_1\to
 (\bE_6\splus\bA_5\splus2\bA_2)\splus\bA_2\splus\bA_1$
(\No$18$);
\dash
$(\bE_6\splus\bA_5\splus2\bA_2)\splus\bA_3\splus\bA_1\to
 (\bE_6\splus4\bA_2)\splus\bA_3\splus\bA_1$
(\No$20'$).
\endDashes
Using equivalence of dessins, see~\cite{degt.kplets}, one can
easily see
that the curves obtained are deformation equivalent to
perturbations of the set of singularities
$\bE_6\splus\bA_5\splus4\bA_2$ (\No$39$ in Table~\ref{tab.e6}).
Their fundamental groups are found in~\cite{degt.8a2}.
Alternatively,\mnote{two sentences added} the second case
(as well as any perturbation of \Nos$12$, $13$, $18$, and~$40$)
is covered by
Corollary~3.10 in~\cite{degt.2a8}, and in the first case one can
show that $\pi_1=\bar\BG3$
similar to \ref{s.loop.torus}. (In the last case,
the group is different:
$\pi_1=\BG4/\Gs_1^2\Gs_2\Gs_3^2\Gs_2$, see~\cite{degt.8a2}, which
agrees with the statement.)

It\mnote{new paragraph added, Remarks 4.5.2 and 4.5.3 removed}
remains to consider perturbations of the distinguished~$\bE_6$
type point~$P$. The only `new' (\ie, not considered
in~\cite{degt.8a2} or~\cite{degt.e6}) set of singularities
realized by irreducible curves of torus type with
$\pi_1\ne\bar\BG3$ is $(\bE_6\splus\bA_5\splus2\bA_2)\splus\bA_4$
(\Nos$9$ and~$41$ in the table), and the `worst' perturbation is
the one in~\iref{s.pert.P}{2a2+a1}, resulting in the extra relation
$[\Ga_1,\Gr^{-2}\Ga_1\Gr^2]=1$. Adding this relation to the
presentations of~$\pi_1$, see~\ref{s.e6+a5+2a2+a4}
and~\iref{s.group.hexagon}{41}, respectively,
and using \GAP~\cite{GAP} similar
to~\ref{s.loop.torus}, one finds that the resulting groups
are isomorphic to~$\bar\BG3$.
(For \No$41$, to make the approach work, one
also needs
to add the relation $(\Ga_1\Ga_2)\1\Ga_2(\Ga_1\Ga_2)=\Ga_3$,
which is known to hold in the group, see
Remark~\ref{rem.extra}.)
\qed

%

\widestnumber\key{EO1}
\refstyle{C}

\enddocument